\documentclass[11pt]{article}
\usepackage{amssymb}
\usepackage{fullpage}
\usepackage{xy}
\xyoption{all}

\title{{\bf Special symplectic connections}}

\date{\mbox{ }}

\def\ftnote#1{\def\footnotemark{}\footnote{#1}\setcounter{footnote}{0}}

\author{Michel Cahen${}^{1,3}$\ftnote{\kern-4.5pt${}^{1}$Universit\'e 
Libre de Bruxelles, Campus Plaine, CP 216, 1050 Bruxelles, Belgium. e-mail: 
mcahen@ulb.ac.be}
\mbox{\hspace{1cm}} Lorenz
J. Schwachh\"ofer${}^{2,3}$\ftnote{\kern-4.5pt${}^{2}$Universit\"at 
Dortmund, Vogelpothsweg 87, 44221 Dortmund, Germany. e-mail: 
lschwach@math.uni-dortmund.de}
\ftnote{\kern-4.5pt${}^{3}$Both authors were supported by the 
Communaut\'e fran\c{c}aise de Belgique, through an Action de Recherche 
Concert\'ee de la Direction de la Recherche Scientifique. The second author 
was also supported through the Schwer\-punkt\-pro\-gramm Globale 
Differentialgeometrie of the Deutsche Forschungsgesellschaft.}}

\begin{document}

\maketitle


\newtheorem{thm}{Theorem}[section]
\newtheorem{lem}[thm]{Lemma}
\newtheorem{prop}[thm]{Proposition}
\newtheorem{df}[thm]{Definition}
\newtheorem{keydf}[thm]{Key Definition}
\newtheorem{cor}[thm]{Corollary}
\newtheorem{rem}[thm]{Remark}
\newtheorem{ex}[thm]{Example}
\newenvironment{proof}{\medskip
\noindent {\bf Proof.}}{\hfill \rule{.5em}{1em}\mbox{}\bigskip}

\def\eps{\varepsilon}

\def\GlR#1{\mbox{\rm Gl}(#1,\R)}
\def\glR#1{{\frak gl}(#1,\R)}
\def\GlC#1{\mbox{\rm Gl}(#1,\C)}
\def\glC#1{{\frak gl}(#1,\C)}
\def\Gl#1{\mbox{\rm Gl}(#1)}

\def\ov{\overline}
\def\ot{\otimes}
\def\und{\underline}
\def\w{\wedge}
\def\ra{\rightarrow}
\def\lra{\longrightarrow}

\newcommand{\SL}{{\mbox{\rm SL}}}
\newcommand{\SO}{{\mbox{\rm SO}}}
\newcommand{\CO}{{\mbox{\rm CO}}}
\newcommand{\GL}{{\mbox{\rm GL}}}
\newcommand{\Sp}{{\mbox{\rm Sp}}}
\newcommand{\U}{{\mbox{\rm U}}}
\newcommand{\SU}{{\mbox{\rm SU}}}
\newcommand{\Spin}{{\mbox{\rm Spin}}}
\newcommand{\E}{{\mbox{\rm E}}}
\newcommand{\G}{{\mbox{\rm G}}}
\newcommand{\T}{{\mbox{\rm T}}}
\newcommand{\Aut}{{\mbox{\rm Aut}}}
\newcommand{\End}{{\mbox{\rm End}}}
\newcommand{\eH}{{\mbox{\rm H}}}
\newcommand{\eP}{{\mbox{\rm P}}}
\newcommand{\K}{{\mbox{\rm K}}}
\renewcommand{\S}{{\mbox{\rm S}}}
\newcommand{\Ad}{{\mbox{\rm Ad}}}
\newcommand{\ad}{{\mbox{\rm ad}}}

\newcommand{\al}{\alpha}
\newcommand{\be}{\beta}
\newcommand{\ga}{\gamma}
\newcommand{\la}{\lambda}
\newcommand{\om}{\omega}
\newcommand{\Om}{\Omega}
\renewcommand{\th}{\theta}
\newcommand{\Th}{\Theta}
\renewcommand{\phi}{\varphi}
\def\rk{\mbox{\rm rk}}

\def\pair#1#2{\left<#1,#2\right>}
\def\big#1{\displaystyle{#1}}

\def\N{{\Bbb N}}
\def\Z{{\Bbb Z}}
\def\R{{\Bbb R}}
\def\C{{\Bbb C}}
\def\bH{{\Bbb H}}
\def\CP{{\Bbb C} {\Bbb P}}
\def\HP{{\Bbb H} {\Bbb P}}
\def\P{{\Bbb P}}
\def\Q{{\Bbb Q}}
\def\F{{\Bbb F}}

\def\O{{\cal O}}
\def\H{{\cal H}}
\def\V{{\cal V}}
\def\W{{\cal W}}
\def\cF{{\cal F}}
\def\cC{{\cal C}}
\def\cR{{\cal R}}

\def\so{{\frak {so}}}
\def\co{{\frak {co}}}
\def\su{{\frak {su}}}
\def\uu{{\frak {u}}}
\def\sl{{\frak {sl}}}
\def\sp{{\frak {sp}}}
\def\e{{\frak {e}}}
\def\csp{{\frak {csp}}}
\def\spin{{\frak {spin}}}
\def\g{{\frak g}}
\def\h{{\frak h}}
\def\k{{\frak k}}
\def\m{{\frak m}}
\def\n{{\frak n}}
\def\t{{\frak t}}
\def\s{{\frak s}}
\def\z{{\frak z}}
\def\p{{\frak p}}
\def\L{{\frak L}}
\renewcommand{\l}{{\frak l}}
\def\X{{\frak X}}
\def\gl{{\frak {gl}}}
\def\hol{{\frak {hol}}}
\renewcommand{\frak}{\mathfrak}
\renewcommand{\Bbb}{\mathbb}

\def\be{\begin{equation}}
\def\ee{\end{equation}}
\def\bi{\begin{enumerate}}
\def\ei{\end{enumerate}}
\def\ba{\begin{array}}
\def\ea{\end{array}}
\def\bea{\begin{eqnarray}}
\def\eea{\end{eqnarray}}
\def\ben{\begin{enumerate}}
\def\een{\end{enumerate}}

\newcommand\A{{\mbox{\bf A}}}
\newcommand\B{{\mbox{\bf B}}}

\renewcommand\b{\mbox{\bf b}}
\renewcommand\c{\mbox{\bf c}}
\def\dq{\slash \!\!\!\! \slash}
\def\codim{\mbox{\rm codim}}
\def\mod{\ \mbox{\rm mod}\ }
\def\hook{\mbox{}\begin{picture}(10,10)\put(1,0){\line(1,0){7}}
  \put(8,0){\line(0,1){7}}\end{picture}\mbox{}}

\begin{abstract}
By a special symplectic connection we mean a torsion free connection which is 
either the Levi-Civita connection of a Bochner-K\"ahler metric of arbitrary 
signature, a Bochner-bi-Lagrangian connection, a connection of Ricci type or 
a connection with special symplectic holonomy. A manifold or orbifold with 
such a connection is called special symplectic.

We show that the symplectic reduction of (an open cell of) a parabolic 
contact manifold by a symmetry vector field is special symplectic in a 
canonical way. Moreover, we show that any special symplectic manifold or 
orbifold is locally equivalent to one of these symplectic reductions.

As a consequence, we are able to prove a number of global properties,
including a classification in the compact simply connected case.
\end{abstract}



\section{Introduction}

Among the basic objects of interest in differential geometry are connections 
on a differentiable manifold $M$ which are compatible with a given geometric 
structure, and the relation between the local invariants of 
such connections and the geometric and topological features of $M$. For 
example, in Riemannian geometry, the Levi-Civita connection of the 
metric is uniquely determined, hence every feature of the connection reflects 
a property of the metric structure.

In contrast, for a symplectic manifold $(M, \om)$, there are many symplectic 
connections, where we call a connection on $M$ symplectic if it is torsion 
free and $\om$ is parallel. Indeed, the space of symplectic connections on 
$M$ is an affine space whose linear part is given by the sections in 
$S^3(TM)$. Thus, in order to investigate 'meaningful' symplectic connections, 
we have to impose further conditions.

In this article, we shall introduce the notion of a {\em special symplectic
connection} which is defined as a symplectic connection on a manifold of
dimension at least $4$ which belongs to one of the following classes.

\bi
\item {\bf Bochner-K\"ahler and Bochner-bi-Lagrangian connections}

If the symplectic form is the K\"ahler form of a (pseudo-)K\"ahler metric, 
then its curvature decomposes into the Ricci curvature and the Bochner 
curvature (\cite{Bo}). If the latter vanishes, then (the Levi-Civita 
connection of) this metric is called Bochner-K\"ahler. 

Similarly, if the manifold is equipped with a bi-Lagrangian structure, i.e. 
two complementary Lagrangian distributions, then the curvature of a 
symplectic connection for which both distributions are parallel decomposes 
into the Ricci curvature and the Bochner curvature. Such a connection is 
called Bochner-bi-Lagrangian if its Bochner curvature vanishes. 

For results on Bochner-K\"ahler and Bochner-bi-Lagrangian connections, see 
\cite{Bochner} and \cite{K} and the references cited therein.

\item {\bf Connections of Ricci type}

Under the action of the symplectic group, the curvature of a symplectic 
connection decomposes into two irreducible summands, namely the Ricci 
curvature and a Ricci flat component. If the latter component vanishes, then 
the connection is said to be of Ricci type. 

Connections of Ricci type are critical points of a certain functional on the 
moduli space of symplectic connections (\cite{Michel}). Furthermore, the 
canonical almost complex structure on the twistor space induced by a 
symplectic connection is integrable iff the connection is of Ricci type 
(\cite{BR}, \cite{V2}). For further properties see also \cite{Michel2}, 
\cite{Michel3}, \cite{Michel4}, \cite{CGS}.

\item {\bf Connections with special symplectic holonomy}

A symplectic connection is said to have {\em special symplectic holonomy} if 
its holonomy is contained in a proper absolutely irreducible subgroup of the 
symplectic group. 

The special symplectic holonomies have been classified in \cite{MS} and 
further investigated in \cite{Br1}, \cite{CMS}, \cite{Habil}, \cite{S}, 
\cite{Advances}.
\ei

We can consider all of these conditions also in the complex case, i.e. for 
complex manifolds of complex dimension at least $4$ with a holomorphic 
symplectic form and a holomorphic connection.

At first, it may seem unmotivated to collect all these structures in
one definition, but we shall provide ample justification for doing so. 
Indeed, our main results show that there is a beautiful link between special 
symplectic connections and parabolic contact geometry. 

For this, consider a (real or complex) simple Lie group $\G$ with Lie algebra 
$\g$. We say that $\g$ is $2$-gradable, if $\g$ contains the root space of a 
long root. In this case, the projectivization of the adjoint orbit of a 
maximal root vector $\cC \subset \P^o(\g)$ carries a canonical $\G$-invariant 
contact structure. Here, $\P^o(V)$ denotes the set of oriented lines through 
$0$ of a vector space $V$, so that $\P^o(V)$ is a sphere if $V$ is real and a 
complex projective space if $V$ is complex. Each $a \in \g$ induces an action 
field $a^*$ on $\cC$ with flow $\T_a := \exp(\F a) \subset \G$, where $\F = 
\R$ or $\C$, which hence preserves the contact structure on $\cC$. Let $\cC_a 
\subset \cC$ be the open subset on which $a^*$ is positively transversal to 
the contact distribution. We can cover $\cC_a$ by open sets $U$ such that the 
local quotient $M_U := \T_a^{loc} \backslash U$, i.e. the quotient of $U$ by 
a sufficiently small neighborhood of the identity in $\T_a$, is a manifold. 
Then $M_U$ inherits a canonical symplectic structure. Our first main result 
is the following

\

\noindent
{\bf Theorem A:} {\em
Let $\g$ be a simple $2$-gradable Lie algebra with $\dim \g \geq 14$, and let 
$\cC \subset \P^o(\g)$ be the projectivization of the adjoint orbit of a 
maximal root vector. Let $a \in \g$ be such that $\cC_a \subset \cC$ is 
nonempty, and let $\T_a = \exp(\F a) \subset \G$. If for an open subset $U 
\subset \cC_a$ the local quotient $M_U = \T_a^{loc} \backslash U$ is a 
manifold, then $M_U$ carries a special symplectic connection.}

\

The dimension restriction on $\g$ guarantees that $\dim M_U \geq 4$ and rules
out the Lie algebras of type $A_1$, $A_2$ and $B_2$.

The type of special symplectic connection on $M_U$ is determined by the Lie
algebra $\g$. In fact, there is a one-to-one correspondence between the 
various conditions for special symplectic connections and simple $2$-gradable 
Lie algebras. More specifically, if the Lie algebra $\g$ is of type $A_n$, 
then the connections in Theorem~A are Bochner-K\"ahler of signature $(p,q)$ if
$\g = \su(p+1,q+1)$ or Bochner-bi-Lagrangian if $\g = \sl(n,\F)$; if $\g$ is 
of type $C_n$, then $\g = \sp(n,\F)$ and these connections are of Ricci type; 
if $\g$ is a $2$-gradable Lie algebra of one of the remaining types, then 
the holonomy of $M_U$ is contained in one of the special symplectic holonomy 
groups. Also, for two elements $a, a' \in \g$ for which $\cC_a, \cC_{a'} 
\subset \cC$ are nonempty, the corresponding connections from Theorem~A are 
equivalent iff $a'$ is $\G$-conjugate to a positive multiple of $a$.

If $\T_a \cong S^1$ then $\T_a \backslash \cC_a$ is an orbifold which carries 
a special symplectic orbifold connection by Theorem~A. Hence it may be 
viewed as the ``standard orbifold model'' for (the adjoint orbit of) $a \in 
\g$. For example, in the case of positive definite Bochner-K\"ahler metrics, 
we have $\cC \cong S^{2n+1}$, and for connections of Ricci type, we have $\cC 
\cong \R\P^{2n+1}$. Thus, in both cases the orbifolds $\T_a \backslash \cC$ 
are weighted projective spaces if $\T_a \cong S^1$, hence the standard 
orbifold models $\T_a \backslash \cC_a \subset \T_a \backslash \cC$ are open 
subsets of weighted projective spaces.

Surprisingly, the connections from Theorem~A exhaust {\em all} special 
symplectic connections, at least locally. Namely we have the following

\

\noindent
{\bf Theorem B:} {\em
Let $(M, \om)$ be a (real or complex) symplectic manifold with a special 
symplectic connection of class $C^4$, and let $\g$ be the Lie algebra associated to the special 
symplectic condition as above.
\bi
\item
Then there is a principal $\hat \T$-bundle $\hat M \ra M$, where $\hat \T$ is 
a one dimensional Lie group which is not necessarily connected, and this 
bundle carries a principal connection with curvature $\om$.
\item
Let $\T \subset \hat \T$ be the identity component.
Then there is an $a \in \g$ such that $\T \cong \T_a \subset \G$, and a
$\T_a$-equivariant local diffeomorphism $\hat \imath: \hat M \ra \cC_a$ which 
for each sufficiently small open subset $V \subset \hat M$ induces a 
connection preserving diffeomorphism $\imath: \T^{loc} \backslash V \ra 
\T_a^{loc} \backslash U = M_U$, where $U := \hat \imath(V) \subset \cC_a$ 
and $M_U$ carries the connection from Theorem~A.
\ei}

The situation in Theorem~B can be illustrated by the following commutative 
diagram, where the vertical maps are quotients by the indicated Lie groups, 
and $\T \backslash \hat M \ra M$ is a regular covering.

\be \label{diagram2}
\xymatrix{
& \hat M \ar[d]^{\mbox{\footnotesize $\T$}} \ar[r]^{\hat \imath}
\ar[dl]_{\mbox{\footnotesize $\hat \T$}} & 
\cC_a \ar[d]^{\mbox{\footnotesize $\T_a$}} 
\\
M & \T \backslash \hat M \ar[r]^{\imath} \ar[l] & \T_a \backslash \cC_a 
}
\ee

In fact, one might be tempted to summarize Theorems~A and B by saying that for
each $a \in \g$, the quotient $\T_a \backslash \cC_a$ carries a canonical 
special symplectic connection, and the map $\imath: \T \backslash \hat M \ra 
\T_a \backslash \cC_a$ is a connection preserving local diffeomorphism. If 
$\T_a \backslash \cC_a$ is a manifold or an orbifold, then this is indeed 
correct. In general, however, $\T_a \backslash \cC_a$ may be neither 
Hausdorff nor locally Euclidean, hence one has to formulate these results 
more carefully. 

As consequences, we obtain the following

\

\noindent
{\bf Corollary C:} {\em
All special symplectic connections of $C^4$-regularity are analytic, and the 
local moduli space of these connections is finite dimensional, in the sense 
that the germ of the connection at one point up to $3$rd order determines the 
connection entirely. In fact, the generic special symplectic connection
associated to the Lie algebra $\g$ depends on $(\rk(\g) - 1)$ parameters.

Moreover, the Lie algebra $\s$ of vector fields on $M$ whose flow preserves 
the connection is isomorphic to $stab(a)/(\F a)$, $\F = \R$ or $\C$, with $a 
\in \g$ from Theorem~B, where $stab(a) = \{ x \in \g \mid [x,a] = 0\}$. 
In particular, $\dim \s \geq \rk(\g) - 1$ with equality implying that $\s$ is 
abelian.
}

\

When counting the parameters in the above corollary, we regard homothetic 
special symplectic connections as equal, i.e. $(M, \om, \nabla)$ is 
considered equivalent to $(M, e^{t_0} \om, \nabla)$ for all $t_0 \in \F$.

We can generalize Theorem~B and Corollary~C easily to orbifolds. Indeed, if 
$M$ is an orbifold with a special symplectic connection, then we can write $M 
= \hat \T \backslash \hat M$ where $\hat M$ is a manifold and $\hat \T$ is a 
one dimensional Lie group acting properly and locally freely on $\hat M$, and 
there is a local diffeomorphism $\hat \imath: \hat M \ra \cC_a$ with the 
properties stated in Theorem~B.

While the analyticity of the connection and the determinedness by the $3$rd
order germ at a point has been known in the Bochner-K\"ahler and 
Bochner-bi-Lagrangian case (\cite{Bochner}\footnote{The $C^4$-regularity of 
the connection is equivalent to the $C^5$-regularity of the Bochner-K\"ahler
metric.}) and for connections with special symplectic holonomies (e.g. 
\cite{CMS}, \cite{MS}), it was unclear what the maximal analytic 
continuations of these structures look like and in which cases they are 
regular. This question is now answered in principle. Furthermore, the 
inequality $\dim \s \geq \rk(\g) - 1$ was known for the Bochner cases 
(\cite{Bochner}), whereas for the special symplectic holonomies, it was only 
known that $\s \neq 0$ (\cite{Advances}).

We also address the question of the existence of compact manifolds with 
special symplectic connections. In the simply connected case, compactness 
already implies that the connection is hermitian symmetric. More 
specifically, we have the following

\

\noindent
{\bf Theorem D:} {\em
Let $M$ be a compact simply connected manifold with a special symplectic
connection of class $C^4$. Then $M$ is equivalent to one of the following 
hermitian symmetric spaces.

\bi
\item
$M \cong (\C\P^p \times \C\P^q, ((q+1) g_0, -(p+1) g_0))$, where $g_0$ is the 
Fubini-Study metric. These are Bochner-K\"ahler metrics of signature $(p,q)$.
Moreover, $M \cong (\CP^n, g_0)$ is also of Ricci type.

\item
$M \cong \SO(n+2)/(\SO(2) \cdot \SO(n))$, whose holonomy is contained in 
the special symplectic holonomy group $\SL(2,\R) \cdot \SO(n) \subset 
\Aut(\R^2 \ot \R^n)$.
\item
$M \cong \SU(2n+2)/\mbox{\rm S}(\U(2) \cdot \U(2n))$, whose holonomy is contained in 
the special symplectic holonomy group $\Sp(1) \cdot \SO(n,\bH) \subset 
\Aut(\bH^n)$.
\item
$M \cong \SO(10)/\U(5)$, whose holonomy is contained in 
the special symplectic holonomy group $\SU(1,5) \subset \GL(20,\R)$.
\item
$M \cong \E_6/(\U(1) \cdot \Spin(10))$, whose holonomy is contained in 
the special symplectic holonomy group $\Spin(2,10) \subset \GL(32, \R)$.
\ei

%
}

\

In particular, there are no compact simply connected manifolds with any of the
remaining types of special symplectic connections, i.e. $M$ can be neither
complex with a holomorphic connection, nor Bochner-bi-Lagrangian, nor can the
holonomy of $M$ be contained in any of the remaining special symplectic
holonomies.

The only case for which Theorem~D was previously known are the positive
definite K\"ahler metrics. In fact, it is shown in \cite{Bochner} that a 
compact positive definite Bochner-K\"ahler manifold must be a quotient of 
$\CP^r \times (\CP^s)^*$, where the asterisk denotes the non-compact dual.

Following this introduction, we first develop the algebraic formulas needed 
to describe the curvature conditions for special symplectic connections 
uniformly. In section~\ref{sec:construction}, we construct the special 
symplectic connections on the local quotients $\T_a^{loc} \backslash \cC_a$ 
and hence prove Theorem~A, and in section~\ref{sec:develope}, we investigate 
the structure equations of special symplectic connections and derive results 
which culminate in Theorem~B. Finally, in the last section we show the 
existence of connection preserving vector fields and Corollary~C, and the 
rigidity result from Theorem~D.

We are grateful to R.Bryant for helpful discussions about Bochner-K\"ahler 
and Bochner-bi-Lagrangian structures, and for valuable comments on the link 
to parabolic contact geometry. Also, it is a pleasure to thank P.Bieliavski, 
S.Gutt and W.Ziller for many stimulating conversations and 
helpful comments. We also thank the referee for many helpful remarks and comments 
which greatly helped to improve this article.

\section{Algebraic preliminaries}

\subsection{A brief review of representation theory}

In this section, we shall give a brief outline of standard facts of 
representation theory of complex semi-simple Lie algebras. For a more 
detailed exposition, see e.g. \cite{FH}, \cite{Hu} or \cite{OV}.

Let $\g_\C$ be a semi-simple complex Lie algebra, and let $\t \subset \g_\C$ 
be a {\em Cartan subalgebra}, i.e. a maximal abelian self-normalizing 
subalgebra. The {\em rank of $\g_\C$} is by definition $\rk(\g_\C) := \dim 
\t$.

If $\rho:\g_\C \ra \End(V)$ is a representation of
$\g_\C$ on a complex vector space $V$, then for any $\la \in \t^*$ we define 
the {\em weight space $V_\la$} by
\[
V_{\la}=\{ v \in V \mid \rho(h)v = \la(h)v\ \mbox{for all}\ h\in \t\}.
\] 
An element $\la \in \t^*$ is called a {\em weight} of $V$ if $V_{\la} \neq
 0$. We let $\Phi \subset \t^*$ be the set of weights of $\rho$, and thus 
have the decomposition
\[
V = \bigoplus_{\la \in \Phi} V_\la.
\]

\noindent
In particular, if $V = \g_\C$ and $\rho$ is the adjoint representation, then 
we get the {\em root decomposition}
\[
\g_\C = \t \oplus \bigoplus_{\al \in \Delta} \g_\al,
\]
i.e. $\t$ is the weight space of weight $0$, and $\Delta \subset \t^*$ is 
the set of non-zero weights. $\Delta$ is called the {\em set of roots} or 
the {\em root system} of $\g_\C$. It is well known that $\dim \g_\al = 1$ for
 all $\al \in \Delta$.
For any root $\al \in \Delta$, there is a unique element $H_{\al} \in
[\g_{\al}, \g_{-\al}] \subset \t$ such that $\al(H_{\al})=2$. 

There is an $\ad(\g_\C)$-invariant non-degenerate symmetric bilinear 
form $B$ on $\g_\C$, the so-called {\em Killing form}, which is given by 
$B(x,y) := \mbox{tr}(\ad_x \circ \ad_y)$ for all $x,y \in \g_\C$. 
We shall use it to identify $\g_\C$ and $\g_\C^*$. The restriction 
of $B$ to $\t$ is non-degenerate as well, and $B(H_\al, H_\al) \in \Z^+$ 
for all $\al \in \Delta$. In fact, there are at most two possible values for
$B(H_\al, H_\al)$ for $\al \in \Delta$ which allows us to speak of {\em long} 
and {\em short} roots, respectively.

Given an element $\la \in \t^*$ and a root $\al$, we let
\be \label{eq:pair}
\pair \la \al := \la(H_\al),\ \ \ \ \ \mbox{so that}\ \ \ \ \ \pair \la \al 
= \frac {2 B(\la, \al)} {B(\al, \al)}.
\ee
Note that $\pair {\ }{\ }$ is linear in the first entry only. We define the 
{\em weight lattice} $\Lambda \subset \t^*$ as the set of elements $\la \in 
\t^*$ such that $\pair \la \al \in \Z$ for all $\al \in \Delta$. Then $\Phi 
\subset \Lambda$ for any representation $\rho$.

For $\la \in \Phi$, the significance of $\pair \la \al \in \Z$ is
the following. If $\la$ occurs as the weight of an irreducible 
representation
of $\g_\C$ and $\pair \la \al > 0$ ($\pair \la \al < 0$, respectively) then 
$\la
- k \al$ ($\la + k \al$, respectively) is also a weight of that
representation for $k = 1, \ldots, |\pair \la \al|$.

For any root $\al \in \Delta$, denote by $\sigma_\al$
the orthogonal reflection of $\t^*$ in the hyperplane perpendicular to 
$\al$. The {\em Weyl group}\, $W$ of $\g_\C$ is the group generated by all
$\sigma_\al$. $W$ is always finite. If $\g_\C$ is simple then $W$ acts
irreducibly on $\t^*$ and transitively on the set of roots of
equal length. The set of weights $\Phi$ of any representation
is $W$-invariant.

If $\g_\C$ is simple, then the adjoint representation $\rho: \g_\C \ra 
\End(\g_\C)$ is irreducible. Also, $|\pair \al 
\beta| \leq 3$ for all $\al, \beta \in \Delta$, and if $\al$ is long and 
$\beta$ short, then either $\pair \al \beta = \pair \beta \al = 0$, 
or $|\pair \al \beta| > 1$ and $|\pair \beta \al| = 1$. Moreover, if $\al$ 
is long then $|\pair \beta \al | \leq 2$, and $\pair \beta \al = \pm 2$ 
iff $\beta = \pm \al$.

\subsection{Special symplectic representations} \label{sec:algebra}

Let $\g_\C$ be a complex simple Lie algebra and let $\G_\C$ be a 
connected complex Lie group with Lie algebra $\g_\C$. Choose a root 
decomposition of $\g$ as in the preceding section, and fix a long root $\al$ 
and an element $0 \neq x \in \g_\al$. Then the orbit of $x$ under the adjoint 
action of $\G_\C$ is called the {\em root cone of $\g_\C$}. 
Evidently, the root cone is well defined, independently of the choice 
of root decomposition. Elements of the root cone are called {\em 
maximal root elements}.

\begin{df} \label{df:parabolic}
Let $\g$ be a simple real or complex Lie algebra. We say that $\g$ is
{\em $2$-gradable} if either $\g$ is complex, or $\g$ is real and contains 
a maximal root element of the simple complex Lie algebra $\g_\C := \g \ot \C$.
\end{df}

We shall justify this terminology in (\ref{eq:decomposeg}) below.
If $\g$ is $2$-gradable and $\G$ is a Lie group with Lie algebra $\g$, then 
we write
\be \label{eq:cones}
\hat \cC := \Ad_G x \subset \g,
\ee
where $x \in \g$ is a maximal root element. Given $x \in \hat \cC$, 
there is a $y \in \hat \cC$ with $B(x,y) \neq 0$, and we can choose 
a root decomposition of $\g$ such that $x \in \g_{\al_0}$ and $y 
\in \g_{-\al_0}$, where $\al_0$ is a long root. Hence $H_{\al_0} \in 
\F [x,y] \subset \t$, so that $\g$ contains the Lie subalgebra 
$\sl_{\al_0} := span<\g_{\al_0}, \g_{-\al_0}, H_{\al_0}>$ which is 
isomorphic to $\sl(2,\F)$, $\F = \R$ or $\C$. Then 
$\ad(H_{\al_0})|_{\g_\beta} = \pair \beta {\al_0} Id_{\g_\beta}$, and since 
$\al_0 \in \Delta$ is a long root, the eigenvalues of $\ad(H_{\al_0})$ are 
$\{0, \pm1, \pm2\}$, so that we get the eigenspace decomposition
\be \label{eq:decomposeg}
\g = \g^{-2} \oplus \g^{-1} \oplus \g^0 \oplus \g^1 \oplus 
\g^2,\ \ \ \ 
\mbox{and}\ \ \ \ [\g^i, \g^j] \subset \g^{i+j},
\ee
where $\g^i = \bigoplus_{\{\beta \in \Delta \mid \pair \beta {\al_0} = i\}} 
\g_\beta$ for $i \neq 0$ and $\g^0 = \t \oplus \bigoplus_{\{\beta \in \Delta 
\mid \pair \beta {\al_0} = 0\}} \g_\beta$. In particular, $\g^{\pm2} = 
\g_{\pm \al_0}$, and $\g^0 = \F H_{\al_0} \oplus \h$, where the Lie algebra 
$\h$ is characterized by $[\h, \sl_{\al_0}] = 0$. Observe that $\g^0$ and 
hence $\h$ are reductive. Thus, as a Lie algebra,
\[
\g^{ev} := \g^{-2} \oplus \g^0 \oplus \g^2 \cong \sl_{\al_0} \oplus 
\h\ \ \ \ \mbox{and}\ \ \ \ \g^{odd} := \g^{-1} \oplus \g^1 
\cong \F^2 \ot V\ \ \ \ \mbox{as a $\g^{ev}$-module},
\]
where $\h$ acts effectively on $V$. Identifying $\h$ with its image 
under this representation, we may regard it as a subalgebra  $\h \subset 
\End(V)$, and hence we have the decomposition
\be \label{eq:symmetricpair}
\g = \g^{ev} \oplus \g^{odd} \cong (\sl(2,\F) \oplus \h) \oplus 
(\F^2 \ot V),
\ee
where this notation indicates the representation $\ad: \g^{ev} \ra 
\End(\g^{odd})$.

We fix a non-zero $\F$-bilinear area form $a \in \Lambda^2 (\F^2)^*$. There 
is a canonical $\sl(2,\F)$-equivariant isomorphim
\be \label{eq:define sl(2)}
S^2(\F^2) \longrightarrow \sl(2,\F),\ \ \ \ (ef) \cdot g := a (e, g) f + a(f, 
g) e\ \ \mbox{for all $e,f,g \in \F^2$},
\ee
and under this isomorphism, the Lie bracket on $\sl(2,\F)$ is given by
\be \label{eq:bracketsl2}
{}[ef, gh] = a(e, g) fh + a(e, h) fg + a(f, g) eh + a(f, h) eg.
\ee
Thus, if we fix a basis $e_+, e_- \in \F^2$ with $a(e_+, e_-) = 1$, then we 
have the identifications
\[
H_{\al_0} = - e_+ e_-,\ \ \ \g^{\pm2} = \F e_\pm^2,\ \ \ \g^{\pm 1} = e_\pm 
\ot V.
\]

\begin{prop} \label{prop:circleproduct}
Let $\g$ be a $2$-gradable simple Lie algebra, and consider the 
decompositions (\ref{eq:decomposeg}) and (\ref{eq:symmetricpair}).
Then there is an $\h$-invariant symplectic form $\om \in \Lambda^2 V^*$ and 
an $\h$-equivariant product $\circ: S^2(V) \ra \h$ such that
\[
{}[\ ,\ ]: \Lambda^2(\g^{odd}) \longrightarrow \g^{ev} \cong \sl(2,\F) 
\oplus \h
\]
is given as
\be \label{eq:bracket}
{}[e \ot x, f \ot y] = \om(x,y) ef + a(e, f) x \circ y\ \ \ \mbox{for 
$e,f \in \F^2$ and $x,y \in V$},
\ee
using the identification $S^2(\F^2) \cong \sl(2,\F) \subset \g^{ev}$ from 
(\ref{eq:define sl(2)}). Moreover, the symmetric bilinear form $(\ ,\ )$ on 
$\g$ defined by
\be \label{eq:Killing}
(u, v) := -\frac 1{2(\dim V + 4)} B(u,v),\ \ \ \ \mbox{for all $u,v 
\in \g$},
\ee
where $B$ is the Killing form of $\g$, satisfies the following:
\bi
\item
$(\g^i, \g^j) = 0$ if $i + j \neq 0$,
\item
$(ef, gh) = a(e, g) a(f, h) + a(e, h) a(f, g)$ for all $e, f, g, h \in 
\F^2$,
\item $B(u,v) = 2\ tr_V(u v) + B_\h(u,v)$ for all $u,v \in \h \subset \g$, 
where $B_\h$ denotes the Killing form of $\h$.
\item
$(e \ot x, f \ot y) = a(e, f) \om(x,y)$, for all $e, f \in \F^2$ and $x, 
y \in V$, using the identification $\g^{odd} \cong \F^2 \ot V$,
\item
For all $x,y,z \in V$ and $h \in \h$, we have
\be \label{eq:Adams}
\ba{c} (h, x \circ y) = \om(hx, y) = \om(hy, x)\\ \\
(x \circ y) z - (x \circ z) y = 2\ \om(y,z) x - \om(x,y) z + \om(x,z) y.
\ea \ee
\ei
\end{prop}

\begin{proof}
By (\ref{eq:decomposeg}) the bracket $[\ ,\ ]: 
\Lambda^2 \g^{odd} \ra \g^{ev}$ is well-defined and must be 
$\g^{ev}$-equivariant by the Jacobi identity. We decompose $\Lambda^2\g^{odd} 
= \Lambda^2(\F^2 \ot V) = S^2(\F^2) \ot \Lambda^2 V \oplus S^2(V)$, so 
that any $\g^{ev}$-equivariant map $\Lambda^2\g^{odd} \ra \g^{ev}$ must be of 
the form (\ref{eq:bracket}) for some $\h$-invariant $\om \in \Lambda^2 V^*$ 
and $\circ: S^2(V) \ra \h$.

Since $(\ ,\ )$ is $\ad_\g$-invariant, i.e. it satisfies the identity 
$([u,v], w) = (u, [v,w])$ for all $u,v,w \in \g$, we have for $u_i \in \g^i$ 
and $u_j \in \g^j$ 
\[
0 = ([H_{\al_0}, u_i], u_j) + (u_i, [H_{\al_0}, u_j]) = (i\ u_i, u_j) + 
(u_i, j\ u_j) = (i+j) (u_i, u_j),
\]
which shows 1. 

To show the second equation, note that the inner product on $S^2(\F^2) \cong 
\sl(2,\F)$ given by the right hand side of this equation is 
$\ad_{\sl(2,\F)}$-invariant and hence must be a multiple of the restriction 
of the Killing form $B$ to $\sl(2,\F)$. Thus, it suffices to verify the 
second equation for $e = g = e_+$ and $f = h = e_-$. In this case, the right 
hand side equals $-1$, whereas the left hand side equals $(e_+ e_-, e_+ e_-) 
= (H_{\al_0}, H_{\al_0})$. But $B(H_{\al_0},H_{\al_0}) = 
tr(\ad(H_{\al_0})^2)$ and since $\ad(H_{\al_0})|_{\g^i} = i Id_{\g^i}$, we
conclude that $(e_+ e_-, e_+ e_-) = -1$ by the choice of the scaling factor in
(\ref{eq:Killing}). This implies 2. Likewise, if $u, v \in 
\h$, then $\ad(u)|_{\sl(2,\F)} = \ad(v)|_{\sl(2,\F)} = 0$, from which 3. follows as well.

For 4. note that $(e_{\pm} \ot x, e_{\pm} \ot y) \in 
(\g^{\pm 1}, \g^{\pm 1}) = 0$ by 1., and from 2. and 
the $\ad_\g$-invariance, we get
\[
(e_+ \ot x, e_- \ot y) = -\frac12 (e_+ \ot x, [e_-^2, e_+ \ot y]) = 
\frac12 ([e_+ \ot x, e_+ \ot y], e_-^2) = \frac12 \om(x,y) \left(e_+^2, 
e_-^2\right) = \om(x,y).
\]
This also implies that $\om$ is symplectic; indeed, if $\om(x, V) = 0$ for 
some $x \in V$, then by 1. and 4. it follows that $(e_+ \ot x, \g) = 0$ so 
that $x = 0$.

To show the first equation in (\ref{eq:Adams}), we note that $(\h, \sl(2,\F)) = 0$ so that 
for $h \in \h$ and $x, y \in V$ we have
\[
(h, x \circ y) = (h, [e_+ \ot x, e_- \ot y]) = ([h, e_+ \ot x], e_- \ot 
y) = (e_+ \ot (hx), e_- \ot y) = \om(hx, y),
\]
where the last identity follows from 4.

Finally, the second equation in (\ref{eq:Adams}) follows when applying the Jacobi identity to 
the elements $e_+ \ot x$, $e_- \ot y$ and $e_- \ot z$.
\end{proof}

In general, given a (real or complex) symplectic vector space $(V, \om)$, 
i.e. $\om \in \Lambda^2 V^*$ is non-degenerate, we define the {\em symplectic 
group} $\Sp(V, \om)$ and the {\em symplectic Lie algebra} $\sp(V, \om)$ by
\[ \ba{c}
\Sp(V, \om) := \{g \in \Aut(V) \mid \om(gx, gy) = \om(x,y) \ \mbox{for all 
$x,y \in V$}\},\\ \\ \sp(V, \om) := \{h \in \End(V) \mid \om(hx, y) + 
\om(x, hy) = 0 \ \mbox{for all $x,y \in V$}\}.
\ea \]
Then $\Sp(V, \om)$ is a Lie group with Lie algebra $\sp(V, \om)$.

\begin{df}
Let $(V, \om)$ be a symplectic vector space over $\F = \R$ or $\C$, and let 
$\h \subset \sp(V, \om)$ be a subalgebra for which there exists an 
$\h$-equivariant map $\circ: S^2(V) \ra \h$ and an $\ad_\h$-invariant inner 
product $(\ ,\ )$ on $\h$ for which the identities (\ref{eq:Adams}) hold. Then we 
call $\h$ a {\em special symplectic subalgebra}. Moreover, we call the 
connected subgroup $\eH \subset \Sp(V, \om)$ with Lie algebra $\h$ a {\em 
special symplectic subgroup}.
\end{df}

Thus, by Proposition~\ref{prop:circleproduct}, each (real or complex)
$2$-gradable simple Lie algebra yields a (real or complex) special symplectic 
subalgebra $\h \subset \End(V)$. The converse is also true. Namely, we have

\begin{prop} \label{prop:special=simple}
Let $(V, \om)$ be a symplectic vector space over $\F = \R$ or $\C$, and 
let $\h \subset \sp(V, \om)$ be a special symplectic subalgebra. Then there 
exists a unique $2$-gradable simple Lie algebra $\g$ over $\F$, which 
admits the decompositions 
(\ref{eq:decomposeg}) and (\ref{eq:symmetricpair}), and the Lie bracket of 
$\g$ is given by (\ref{eq:bracket}).
\end{prop}

\begin{proof}
Given the special symplectic Lie algebra $\h \subset \sp(V, \om)$, we 
define the $(\sl(2, \F) \oplus \h)$-equivariant map $R: \Lambda^2 (\F^2 \ot 
V) \ra \sl(2,\F) \oplus \h$ by (\ref{eq:bracket}) and verify that $R$ 
satisfies the Jacobi identity by the property of $\circ$.

Thus, $R$ defines a Lie algebra structure on $\g := \sl(2, \F) \oplus \h 
\oplus \F^2 \ot V$ which makes $(\g, \sl(2,\F) \oplus \h)$ into a symmetric 
pair. Choose a basis $e_\pm$ of $\F^2$ with $a(e_+, e_-) = 1$ and let 
$\g^0 := \F e_+ e_- \oplus \h$, $\g^{\pm 1} := e_\pm \ot V$ and $\g^{\pm 2} 
:= \F e_\pm^2$. Then $[\g^i, \g^j] \subset \g^{i+j}$ follows from the 
definition of the bracket, so that (\ref{eq:decomposeg}) and 
(\ref{eq:symmetricpair}) hold.

Let $\g' \subset \g$ be an ideal. Since $e_+ e_-$ is a grading element, it 
follows that $\g' = \bigoplus_{i=-2}^2 (\g' \cap \g^i)$. Moreover, $\g' 
\cap \sl(2,\F) \subset \sl(2,\F)$ is an ideal, hence either $\g' \cap 
\sl(2,\F) = 0$ or $\sl(2,\F) \subset \g'$.

First, suppose that $\g' \cap \sl(2,\F) = 0$ so that $\g' \cap \g^{\pm2} = 
0$. If $e_\pm \ot x \in \g' \cap \g^{\pm 1}$, then for all $y \in V$, we have 
$[e_\pm \ot x, e_\pm \ot y] = \om(x,y) e_\pm^2 \in \g' \cap \g^{\pm 2} = 0$ 
so that $\om(x,y) = 0$ for all $y \in V$, i.e. $x = 0$, hence $\g' \cap 
\g^{\pm 1} = 0$, whence $\g' \subset \g_0$. Next, $[\g', \g^{\pm 2}] 
\subset \g' \cap [\g^0, \g^{\pm 2}] = \g' \cap \g^{\pm 2} = 0$, so that 
$\g' \subset \h$. Finally, if $h \in \g' \subset \h$, then for all $x \in V$, 
$[h, e_\pm \ot x] = e_\pm \ot (hx) \in \g' \cap \g^{\pm 1} = 0$, i.e. $hx = 
0$ for all $x \in V$, hence $h = 0$, i.e. $\g' = 0$.

On the other hand, if $\sl(2,\F) \subset \g'$, then $e_+ e_- \in \g'$ so 
that $\g^i = [e_+ e_-, \g^i] \subset \g'$ for all $i \neq 0$. Moreover, 
$[\g^1, \g^{-1}] \subset \g'$, so that $x \circ y \in \g'$ for all $x,y \in 
V$. By the first identity of (\ref{eq:Adams}), we have $V \circ V = \h$, so 
that $\h \subset \g'$ and hence $\g' = \g$.

We conclude that $\g$ is simple, and since $\ad(e_+ e_-)$ is diagonalizable, 
we can choose the root subalgebra $\t$ such that $e_+ e_- \in \t$. Then $\t 
= \F e_+ e_- \oplus (\t \cap \h)$, and hence $[\t, \g^{\pm2}] = \g^{\pm2}$, 
so that $\g^{\pm2} = \g_{\pm \al_0}$ are root spaces and $H_{\al_0} = -e_+ 
e_-$. Recall that $\ad(H_{\al_0})|_{\g_\beta} = \pair \beta {\al_0} 
Id_{\g_\beta}$ which implies that $|\pair \beta {\al_0}| \leq 1$ for all 
roots $\beta \neq \pm \al_0$, hence $\al_0$ is a long root.
\end{proof}

{}From this proposition, we obtain a complete classification of special 
symplectic subalgebras by considering all complex simple Lie algebras and 
their $2$-gradable real forms (\cite{OV}).

\begin{cor} \label{cor:sympl}
Table~1 yields the complete list of special symplectic subgroups $\eH 
\subset \Sp(V, \om)$.
\end{cor}

\begin{table}
\begin{center}
{\footnotesize
{\bf Table~1: Special symplectic subgroups}\\
{Notation: $\F = \R$ or $\C$.}
\begin{tabular}{|r|c|c|c|c|c}
\hline
& Type of $\Delta$ & $\G$ & $\eH$ & $V$\\
\hline 
\hline
&&&& \vspace{-3mm}\\
(i) & $\mbox{A}_{k}$, $k \geq 2$ & $\SL(n+2, \F)$, $n \geq 1$ & $\GL(n,\F)$ & 
$W \oplus W^*$ with $W \cong \F^n$\\
&&&& \vspace{-3mm}\\
(ii) & & $\SU(p+1, q+1)$, $p + q \geq 1$ & $\U(p,q)$ & $\C^{p+q}$\\
\hline
&&&& \vspace{-3mm}\\
(iii) & $\mbox{C}_k$, $k \geq 2$ & $\Sp(n+1, \F)$ & $\Sp(n, \F)$ & $\F^{2n}$\\
\hline
&&&& \vspace{-3mm}\\
(iv) & $\mbox{B}_{k}$, $\mbox{D}_{k+1}$, $k \geq 3$ & $\SO(n+4,\C)$, $n\geq 
3$ & $\SL(2,\C) \cdot \SO(n,\C)$ & $\C^2 \ot \C^{n}$\\
&&&& \vspace{-3mm}  \\
(v) & & $\SO(p+2, q+2)$, $p+q \geq 3$ & $\SL(2,\R) \cdot \SO(p,q)$ & $\R^2 
\ot \R^{p+q}$\\
&&&& \vspace{-3mm}\\
(vi) & & $\SO(n+2, \bH)$, $n\geq 2$ & $\Sp(1) \cdot \SO(n,\bH)$ & $\bH^n$\\
\hline
&&&& \vspace{-3mm}\\
(vii) & $\mbox{G}_2$ & $\mbox{G}_2'$, $\mbox{G}_2^\C$ & $\SL(2, \F)$ & 
$S^3(\F^2)$\\
\hline
&&&& \vspace{-3mm}\\
(viii) & $\mbox{F}_4$ & $\mbox{F}_4^{(1)}$, $\mbox{F}_4^\C$ & $\Sp(3,\F)$ & 
$\F^{14} \subset \Lambda^3 \F^6$\\
\hline
&&&& \vspace{-3mm}\\
(ix) & $\mbox{E}_6$ & $\mbox{E}_6^\F$ & $\SL(6,\F)$ &  $\Lambda^3 \F^6$\\
&&&& \vspace{-3mm}\\
(x) & & $\mbox{E}_6^{(2)}$ & $\SU(1,5)$ & $\R^{20} \subset \Lambda^3 \C^6$\\
&&&& \vspace{-3mm}\\
(xi) & & $\mbox{E}_6^{(3)}$ & $\SU(3,3)$ & $\R^{20} \subset \Lambda^3 \C^6$\\
\hline
&&&& \vspace{-3mm}\\
(xii) & $\E_7$ & $\E_7^\C$ & $\Spin(12,\C)$ & $\Delta^\C \cong \C^{32}$\\ 
&&&& \vspace{-3mm}\\
(xiii) & & $\mbox{E}_7^{(5)}$ & $\Spin(6,6) $ & $\R^{32} \subset \Delta^\C$\\ 
&&&& \vspace{-3mm}\\
(xiv) & & $\mbox{E}_7^{(6)}$ & $\Spin(6,\bH) $ & $\R^{32} \subset \Delta^\C$\\
&&&& \vspace{-3mm}\\
(xv) & & $\mbox{E}_7^{(7)}$ & $\Spin(2,10) $ & $\R^{32} \subset \Delta^\C$\\ 
\hline
&&&& \vspace{-3mm}\\
(xvi) & $\mbox{E}_8$ & $\E_8^\C$ & $\E_7^{\C}$ & $\C^{56}$\\
&&&& \vspace{-3mm}\\ 
(xvii) & & $\mbox{E}_8^{(8)}$ & $\E_7^{(5)}$ & $\R^{56}$\\
&&&& \vspace{-3mm}\\
(xviii) & & $\mbox{E}_8^{(9)}$ & $\E_7^{(7)}$ & $\R^{56}$\\
\hline
\end{tabular}
}
\end{center}
\end{table}

It is worth pointing out that in the case $\h = \sp(V, \om)$ the map $\circ:
S^2(V) \ra \h$ is an isomorphism which is given explicitly by
\be \label{eq:circle-spV}
(x \circ y) z = \om(x,z) y + \om(y,z) x\ \ \ \mbox{for all $x,y,z \in V$.}
\ee
Namely, by Proposition~\ref{prop:special=simple} it suffices to show that 
this product is well defined, $\h$-equivariant and satisfies 
(\ref{eq:Adams}), and all of this is easily verified.

\begin{df} \label{df:associated}
Let $\h \subset \sp(V, \om)$ be a special symplectic Lie algebra, and let 
$\g$ be the (unique) simple Lie algebra from 
Proposition~\ref{prop:special=simple}. Then we say that $\h$ is {\em 
associated to $\g$}. Let $\G$ be a connected Lie group with Lie algebra $\g$. 
Then we say that the special symplectic group $\eH \subset \Sp(V, \om)$ is 
{\em associated to $\G$}.
\end{df}

\begin{prop} \label{prop:Hconnected}
Let $\h \subset \sp(V, \om)$ be a special symplectic Lie algebra and 
$\eH \subset \Sp(V, \om)$ be the corresponding special symplectic Lie subgroup. Then 
$\eH \subset \Sp(V, \om)$ is closed and reductive, and
\be \label{eq:circleinv}
\h = \{ h \in \sp(V, \om) \mid [h, x \circ y] = (hx) \circ y + x \circ 
(hy)\ \mbox{for all $x,y \in V$} \}.
\ee

Moreover, let $\g \cong \sl(2,\F) \oplus \h \oplus \F^2 \ot V$ be the simple 
Lie algebra from Proposition~\ref{prop:special=simple} and $\G$ the 
corresponding simply connected Lie group from Definition~\ref{df:associated}. 
Then the Lie subgroup
\be \label{eq:definetildeH}
\tilde \eH := \{g \in G \mid \Ad_g|_{\g^{-2} \oplus \g^2} = Id_{\g^{-2} 
\oplus \g^2}\} \subset G
\ee
is generated by $\eH$ and the center $Z(\G)$.
\end{prop}

\begin{proof}
In principle, we could prove this theorem from Table~1, but we prefer 
to give more conceptual arguments.

Let us suppose that $\h$ and $V$ are complex. Then, by 
Proposition~\ref{prop:special=simple}, we can find a complex simple Lie 
algebra $\g$ for which (\ref{eq:decomposeg}) holds. Thus, $\g^0 = \t 
\oplus \bigoplus_{\{\beta \in \Delta \mid \pair \beta {\al_0} = 0\}} 
g_\beta$ where $\Delta$ is the set of roots of $\g$. Then $\g^0$ is evidently 
reductive, and since $\g^0 \cong \C \oplus \h$, it follows that $\h$ is 
reductive as well, hence so is every real form of $\h$. Thus, $\eH$ 
is also reductive.

Let $\tilde \h$ denote the right hand side of (\ref{eq:circleinv}). 
Then the $\h$-equivariance of $\circ$ implies that $\h \subset 
\tilde \h$. Also, $\h = V \circ V$ by the first identity of 
(\ref{eq:Adams}) so that $\h$ is an ideal of $\tilde \h$. Therefore, if 
$\tilde h \in \tilde \h$ then we define 
\[
\phi: \g \ra \g\ \ \ \ \mbox{ by }\ \ \ \ \phi(\sl(2, \F)) = 0,\ \ \ 
\phi|_\h = (\ad_{\tilde h})|_\h,\ \ \ \phi|_{\F^2 \ot V} := Id_{\F^2} \ot 
\tilde h.
\]
Since $\ad_{\tilde h}(\h) \subset \h$, this definition makes sense. 
Moreover, it is now straightforward to verify that $\phi$ is a derivation of 
$\g$, and since $\g$ is simple, it follows that $\phi = \ad_h$ for some $h 
\in \g$. But $\phi(\sl(2,\F)) = 0$, so that $h \in \h$, hence $e \ot (h x) = 
\ad_h (e \ot x) = \phi(e \ot x) = e \ot (\tilde h x)$ for all $e \in \F^2$ 
and $x \in V$, whence $\tilde h = h \in \h$ which shows (\ref{eq:circleinv}).

Now the subgroup $\{ h \in \Sp(V, \om) \mid \Ad_h (x \circ y) = (hx) \circ 
(hy) \} \subset \Sp(V, \om)$ is closed and has $\h$ as its Lie algebra by 
(\ref{eq:circleinv}), thus $\eH$ is its identity component and hence 
also closed.

For the last part, note that the Lie algebra of $\tilde \eH$ equals
$\{x \in \g \mid [x, \g^{\pm 2}] = 0\} = \h$. As $\eH$ is connected, this 
implies that $\eH \subset \tilde \eH$ is the identity component, and it thus 
suffices to show that every component of $\tilde \eH$ contains an element of 
$Z(\G)$.

Let $g \in \tilde \eH$. Then $\h$ is $\Ad_g$-invariant, and if we let 
$\t_\h \subset \h$ be a Cartan subalgebra of $\h$, so that $\t_\g := \t_\h 
\oplus \F e_+ e_- \subset \g^0$ is a Cartan subalgebra of $\g$, then 
$\Ad_g(\t_\h) \subset \h$ is another Cartan subalgebra. Since any two Cartan 
subalgebras are conjugate via an element of $\eH$, we may assume w.l.o.g. 
that $\Ad_g(\t_\h) = \t_\h$, and since $\Ad_g(e_+ e_-) = e_+ e_-$, it follows 
that $\Ad_g \in Norm(\t_\g)$. Thus, $\Ad_g$ yields an inner automorphism 
of the root system of $\g$ which stabilizes the root $\al_0$, so that the
restriction $(\Ad_g)|_{\t_\h}$ is an inner automorphism of the root system of 
$\h$, hence after multiplying $g$ by an element of $Norm(\t_\h) \subset \eH$, 
we may assume that $(\Ad_g)|_{\t_\g} = Id_{\t_\g}$, so that $g \in \T = 
\exp(\t_\g) = \exp(\F e_+ e_-) \exp(\t_\h)$. Since $\exp(\t_\h) \subset \eH$, 
we may further assume that $g = \exp(t e_+ e_-)$ for some $t \in \F$, hence 
$\Ad_g|_{\g^i} = c^i Id_{\g^i}$ with $c := \exp(-t)$. But $g \in \tilde \eH$, 
so that we must have $c = \pm 1$. 

If $c = 1$ then $\Ad_g = Id$, i.e. $g \in Z(\G)$, so that we are done.

If $c = -1$ then $\F = \C$ and $\Ad_g|_{\g^{\pm 1}} = -Id_{\g^{\pm 1}}$, 
hence we are done if we can show that $-Id_V \in \eH$, since then $g \cdot 
(-Id_V) \in Z(\G)$.

If $\eH = \Sp(V, \om)$, then this is certainly the case, and if $\eH 
\subsetneq \Sp(V, \om)$ is a proper subgroup, then we shall see in 
Lemma~\ref{lem:weights}, 5. that there is an $h \in \t_\h$ such that 
$\la(h)$ is an odd integer for all weights $\la$ of $V$, hence 
$\exp(\sqrt{-1} \pi h) = -Id_V \in \eH$.
\end{proof}

In general, for a given Lie subalgebra $\h \subset \End(V)$ we define the 
space of {\em formal curvature maps} as
\[
K(\h) := \{ R \in \Lambda^2 V^* \ot \h \mid R(x,y) z + R(y,z) x + R(z,x) 
y = 0\ \mbox{for all $x,y,z \in V$}\}.
\]
This terminology is due to the fact that the curvature map of a torsion 
free connection always satisfies the first Bianchi identity, i.e. is 
contained in $K(\h)$ for an appropriate $\h$. $K(\h)$ is an $\eH$-module in 
an obvious way.

There is a map $Ric: K(\h) \ra V^* \ot V^*$, given by $Ric(R)(x,y) := 
tr(R(x, \und\ )y)$ for all $R \in K(\h)$ and $x, y \in V$. Note that 
$Ric(R)(x,y) - Ric(R)(y,x) = tr R(x,y)$. Thus, if $\h \subset \sl(n,\F)$, 
then $Ric(R) \in S^2(V^*)$.

\begin{prop} \label{prop:Rh-injective}
Let $\h \subset \sp(V, \om)$ be a special symplectic subalgebra. 
Then there is an $\eH$-equi\-va\-ri\-ant injective map $\h \ra K(\h)$, given 
by
\be \label{eq:defRh}
h \longmapsto R_h,\ \ \ \mbox{where}\ \ \ R_h(x,y) := 2\ 
\om(x,y) h + x \circ (hy) - y \circ (hx).
\ee
In fact, $Ric(R_h) = 0$ iff $h = 0$.
\end{prop}

\begin{proof}
The fact that $R_h \in K(\h)$ follows immediately from (\ref{eq:Adams}), and 
the $\eH$-equivariance is evident. The injectivity will follow from the last 
statement. We begin with the following two lemmata.

\begin{lem} \label{lem:determineConstant}
Let $\g$ be a (real or complex) semi-simple Lie algebra and let $\h \subset \g$ 
be simple. Then there is a $c \in (0, 1]$ such that for the Killing forms of $\g$ and 
$\h$, the relation
\[
B_\h = c (B_\g)|_\h
\]
holds. Moreover, $c = 1$ iff $\h \lhd \g$.
\end{lem}

\begin{proof} Since $\h$ is simple and both $B_\h$ and $(B_\g)|_\h$ are 
$ad_\h$-invariant, Schur's Lemma implies that this relation holds 
for some $0 \neq c \in \F$. Note that $c$ remains unchanged if we replace 
$\h$ and $\g$ by their complexification or a real form. Thus, it suffices to show  
that $c \in (0, 1]$ for compact Lie algebras $\h \subset \g$, i.e., for $B_\h, 
B_\g < 0$.

For $0 \neq x \in \h$, we have $B_\g(x, x) = trace (ad_x^2) = B_\h(x, x)  + 
trace (ad_x^2|_{\h^\perp})$. Since $ad_x$ is skew symmetric w.r.t. the 
(positive definite) inner product $-B_\g$, it follows that $ad_x^2|_{\h^\perp}$ is 
negative semidefinite, so that $trace (ad_x^2|_{\h^\perp}) \leq 0$ with equality 
iff $ad_x|_{\h^\perp} = 0$, which implies the claim. 
\end{proof}

\begin{lem} \label{lem:Ricci}
Let $\h \subset \sp(V, \om)$ be a symplectic subalgebra. Then 
$Ric(R)(x,y) = -\om(R(\om^{-1})x,y)$. In particular, $Ric(R) \in \h \subset \sp(V) \cong S^2(V)$.
\end{lem}

\begin{proof}
Let $(e_i, f_i)$ be a basis of $V$ such that, using the summation 
convention, $\om^{-1} = e_i \w f_i$. Thus,
\[
\ba{llll}
Ric(R)(x,y) & = & tr(R(x, \und\ )y) = \om(R(x, e_i )y, f_i) - \om(R(x, 
f_i)y, e_i)\\ \\
& = & \om(R(x, e_i )f_i, y) + \om(R(f_i, x)e_i, y) = -\om(R(e_i, f_i) x, y).
\ea
\]

\vspace{-9mm}
\hfill \end{proof} 

\noindent
Let us now suppose that $Ric(R_h) = 0$. By the lemma, this is the case iff 
for all $u \in \h$ we have
\[
\ba{lll}
0 = (R_h(e_i, f_i), u) & = & 2 \om(e_i, f_i) (h, u) + (e_i \circ (h f_i), u) 
- (f_i \circ (h e_i), u)\\ \\
& = & \dim V (h, u) + \om((u h f_i), e_i) - \om((u h e_i), f_i)\\ \\
& = & \dim V (h, u) - tr_V (uh)\\ \\
& = & \dim V (h, u) - \frac12 (B(h, u) - B_\h(h, u))\\ \\
& = & \dim V (h, u) - \frac12 (-2(\dim V + 4) (h,u) - B_\h(h,u))\\ \\
& = & 2 (\dim V + 2) (h,u) + \frac12 B_\h(h,u).
\ea
\]
Here, we use repeatedly the identities from 
Proposition~\ref{prop:circleproduct}.
Let
\[
\h = \h_0 \oplus \h_1 \oplus \ldots \oplus \h_k
\]
be the decomposition of $\h$ with $\h_0 := \z(\h)$ and $\h_r$ simple for $r 
\geq 1$. By simplicity of $\h_r$, there are constants $c_r \in [0,1]$ such 
that $B_{\h_r} = c_r B|_{\h_r}$, where $c_0 = 0$ and $c_r > 0$ for $r > 0$. 
Thus, if we decompose $h = h_0 + \ldots + h_k$ with $h_r \in \h_r$, then $R_h 
= 0$ iff for all $u_r \in \h_r$ we have
\[
\ba{lll}
0 & = & 2 (\dim V + 2) (h_r, u_r) + \frac12 B_{\h_r} (h_r, u_r) = 2 (\dim V + 
2) (h_r, u_r) + \frac12 c_r B(h_r, u_r)\\ \\
& = & (2 (\dim V + 2) - c_r (\dim V + 4)) (h_r, u_r),
\ea
\]
using again Proposition~\ref{prop:circleproduct}. But since $c_r \leq 1$ by Lemma~\ref{lem:determineConstant}, it 
follows that $2 (\dim V + 2) - c_r (\dim V + 4) \geq \dim V > 0$, so that 
we must have $h_r = 0$ for all $r$ which completes the proof.
\end{proof}

For a special symplectic subalgebra $\h \subset \sp(V,\om)$, we can now 
decompose its curvature space as an $\h$-module into
\be \label{eq:defRhspace}
K(\h) = \cR_\h \oplus \W_\h,\ \ \ \ \mbox{where}\ \ \ \ \cR_\h  = \{ R_h 
\mid h \in \h \}.
\ee

By Proposition~\ref{prop:Rh-injective} and Lemma~\ref{lem:Ricci}, it follows 
that $\cR_\h \cong \h$ as an $\eH$-module and $\W_\h$ is the kernel of the 
map $Ric: K(\h) \ra \h \subset \sp(V, \om) \cong S^2(V^*)$, i.e. $\W_\h$ 
consists of all {\em Ricci flat curvature maps}.

In fact, the curvature spaces $K(\h)$ have been calculated. 
Summarizing, we have the following

\begin{thm} \label{thm:curvatures}
Let $\eH \subset \Sp(V, \om)$ be a special symplectic subgroup with Lie 
algebra $\h \subset \sp(V, \om)$ listed in Table~1. Then
\bi
\item
For the representations corresponding to $(i)$ and $(ii)$, we have $\W_\h 
= 0$ if $n = 1$ ($p+q = 1$, respectively) and $\W_\h \neq 0$ if $n \geq 2$ 
($p+q \geq 2$, respectively).
\item
For the representations corresponding to $(iii)$, we have $\W_\h = 0$ for 
$n = 1$ whereas $\W_\h \neq 0$ for $n \geq 2$.
\item
For the representations corresponding to entries $(iv)$ -- $(xviii)$, we have 
$K(\h) = \cR_\h$ and hence $\W_\h = 0$.
\ei
\end{thm}

\begin{proof}
First of all, note that since $\h_\C = \h_\R \ot \C$ and $V_\C = V_\R \ot 
\C$, we also have $K(\h_\C) = K(\h_\R) \ot \C$ and $\cR_{\h_\C} = 
\cR_{\h_\R} \ot \C$ by complexification. Thus, 
it suffices to show the claim for the complex representations.

Therefore, to show the first part, it suffices to show that in case $(i)$, 
$K(\h) \cong S^2(W) \ot S^2(W^*)$ as an $\h$-module, so that the assertion 
follows by a dimension count. To see this, let $x, y \in W$ and $\ov z, 
\ov w \in W^*$. Then for any $R \in K(\h)$ we have $R(\ov z, x) y - R(\ov z, 
y) x = - R(x, y) \ov z$, and since the left hand side lies in $W$ while 
the right hand side lies in $W^*$, it follows that both sides vanish.

The vanishing of the right hand side implies that $R(W, W) = 0$ since $x, y 
\in W$ and $\ov z \in W^*$ are arbitrary. Analogously, $R(W^*, W^*) = 0$. 
Moreover, the vanishing of the left hand side implies that $R(\ov z, x) 
y = R(\ov z, y)x$ and, analogously, $R(x, \ov z) \ov w = R(x, \ov w) \ov z$.
Thus, if we define the tensor $\sigma_R \in W \ot W \ot W^* \ot W^*$ by
\be \label{eq:sigma-R}
\sigma_R(x, y, \ov z, \ov w) := \ov w (R(\ov z, x) y) = - (R(\ov z, x) \ov w) 
y\ \ \ \ \mbox{for all $x, y \in W$ and $\ov z, \ov w \in W^*$},
\ee
then $\sigma_R$ is symmetric in $x$ and $y$ and in $\ov z$ and $\ov w$, i.e. 
$\sigma_R \in S^2(W) \ot S^2(W^*)$.

Conversely, given $\sigma \in S^2(W) \ot S^2(W^*)$, one verifies that the 
map $R_\sigma: \Lambda^2(V) \ra \h$ determined by $R(W, W) = R(W^*, W^*) = 
0$ and (\ref{eq:sigma-R}) lies in $K(\h)$, showing the above equivalence.

For the second part, consider the Koszul exact sequence $\ldots \ra \Lambda^k
V^* \ot S^l(V^*) \ra \Lambda^{k+1} \ot S^{l-1}(V^*) \ra \ldots$ where the maps
are given by skew symmetrization. One observes that under the identification 
$\sp(V, \om) \cong S^2(V^*)$ we may regard $K(\sp(V))$ as the kernel of the 
map $\Lambda^2 V^* \ot S^2 V^* \ra \Lambda^3 V^* \ot V^*$, hence $K(\sp(V))
\cong (V^* \ot S^3(V^*))/S^4(V^*)$, so that the statement follows by a 
dimension count (cf. \cite{Michel}). The last part was shown in \cite{MS}.
\end{proof}

Now the {\em second Bianchi identity} of the covariant derivative of a 
torsion free connection motivates the following definition. We define the 
space of {\em covariant $\cR$-derivations} by
\be \label{eq:defR1}
\cR_\h^{(1)} := \left\{ \psi \in V^* \ot \cR_\h \mid \psi(x)(y,z) + 
\psi(y)(z,x) + \psi(z)(x,y) = 0 \ \ \mbox{for all $x, y, z \in V$} \right\}.
\ee
Again, $\cR_\h^{(1)}$ is an $\eH$-module in an obvious way.

\begin{prop} \label{prop:R1h}
Let $\h \subset \sp(V, \om)$ be a special symplectic subalgebra other than 
the subalgebra $\h = \sl(2,\F)$, $V = \F^2$. Then as an $\h$-module, 
$\cR_\h^{(1)} \cong V$ with an explicit isomorphism given by
\[
u \longmapsto \psi_u,\ \ \ \mbox{where}\ \ \ \psi_u(x) := R_{u \circ x} \in 
\cR_\h\ \ \ \mbox{for all $u, x \in V$.}
\]
\end{prop}

\begin{proof} As in the proof of Theorem~\ref{thm:curvatures}, it suffices 
to show the proposition in the complex case by complexifying $\h$ and $V$.

Using (\ref{eq:Adams}), it is straightforward to verify that 
$\psi_u \in \cR_\h^{(1)}$ for all $u \in V$. Also, $\psi_u = 0$ iff $R_{u 
\circ V} = 0$ iff $u \circ V = 0$ by Proposition~\ref{prop:Rh-injective}. 
But, again by (\ref{eq:Adams}), $u \circ V = 0$ iff $u = 0$, so that 
$\{\psi_u \mid u \in V\} \subset \cR_\h^{(1)}$ is isomorphic to $V$ as an 
$\eH$-module.

If $\h = \sp(V, \om)$ then $\circ: S^2(V) \ra \h$ is given in
(\ref{eq:circle-spV}), and from there the statement follows for $\dim V > 2$ 
by a direct calculation (\cite{Michel}). On the other hand, if $\dim V = 2$ 
then evidently, $\cR^{(1)}_\h = V \ot \h$, and $\dim \h \in \{1,3\}$ as $\h 
\subset \sl(2,\C)$. Thus, by a dimension count the statement follows if $\dim 
\h = 1$ while it fails if $\dim \h = 3$, i.e. if $\h = \sl(2,\C)$ and $V = 
\C^2$.

Thus, the major part of the proof is to show that the inclusion $\{\psi_u \mid 
u \in V\} \subset \cR_\h^{(1)}$ is an equality if $\h \subsetneq \sp(V, \om)$ and 
$\dim V > 2$. For this, we begin with the following

\begin{lem} \label{lem:weights} (cf. \cite{S})
Let $\h \subsetneq \sp(V, \om)$ be a special symplectic proper subalgebra, 
where $\h$ and $V$ are complex and $\dim V > 2$. Let $\t_\h \subset \h$ be a 
Cartan subalgebra and $\Delta_\h$ be the set of roots of $\h$. 
Consider the decomposition $V = \bigoplus_{\la \in \Phi} V_\la$ where $\Phi 
\subset \t_\h^*$ is the set of weights. Then the following holds:
\bi
\item
All weight spaces $V_\la$ are one dimensional, and if $\la \in \Phi$ then 
$-\la \in \Phi$.
\item
There are at most two possible length for the weights which allows to 
refer to {\em long} and {\em short} weights.
\item
If $\la_0 \in \Phi$ is a long weight, then there is a disjoint decomposition
\[
\Phi = \Phi_{-3} \cup \Phi_{-1} \cup \Phi_{1} \cup \Phi_{3},
\ \ \ \ \mbox{where $\Phi_{\pm3} = \{\pm \la_0 \}$ and $\Phi_{\pm 1} = \{ \mu 
\in \Phi \mid \pm \la_0 - \mu \in \Delta_\h\}$}.
\]
\item
Let $V_{\frac i2} := \bigoplus_{\la \in \Phi_i} V_\la$ for $i \in \{\pm 1, 
\pm 3\}$. Then there are decompositions
\[ \ba{c}
\h = \h_{-1} \oplus \h_0 \oplus \h_1,\ \ \ \ V = V_{-\frac32} \oplus 
V_{-\frac12} \oplus V_{\frac12} \oplus V_{\frac32}\ \ \ \ \mbox{with}\\ \\ 
{}[\h_i, \h_j] \subset \h_{i+j},\ \ \ \ \h_i V_r \subset V_{i+r},\ \ \ \ 
V_r \circ V_s \subset \h_{r+s},\ \ \ \ \h_i = \bigoplus_{r+s=i} V_r \circ V_s.
\ea
\]
\item
Let $v_\pm \in V_{\pm \frac32}$, $w_r \in V_r$ and $h_i \in \h_i$. Then
$
(v_+ \circ v_-) w_r = -2r\ \om({v_+},{v_-}) w_r
$ and $[v_+ \circ v_-, h_i] = -2i\ \om(v_+, v_-) h_i$.
\ei
\end{lem}

\begin{proof}
Let $\g$ be the simple Lie algebra associated to $\h$ by 
Proposition~\ref{prop:special=simple}, and let $\Delta$ be the root system of 
$\g$. Note that $\t_\h = \t \cap (H_{\al_0})^\perp$ where $\t$ is the Cartan 
subalgebra of $\g$. Moreover, $\Delta_\h = \{\beta \in \Delta \mid 
\pair \beta {\al_0} = 0\} \subset \Delta$, and $V \cong \g^1 = 
\bigoplus_{\{\beta \in \Delta \mid \pair \beta {\al_0} = 1\}} \g_\beta$ as an 
$\h$-module. It follows that
\[
\Phi = \left\{ \left. \la = \beta - \frac12 \al_0 \right| \beta \in \Delta, 
\pair \beta {\al_0} = 1 \right\}\ \ \ \ \mbox{and}\ \ \ \ V_\la = \g_\beta.
\]
Thus, $\dim V_\la = 1$ as all root spaces are one dimensional. Moreover, 
if $\pair \beta {\al_0} = 1$, then $\gamma := \al_0 - \beta \in \Delta$ 
and $\pair \gamma {\al_0} = 1$, whence $-\la = -(\beta - \frac12 \al_0) =
\gamma - \frac12 \al_0 \in \Phi$.

Next, $(\la, \la) = (\beta - \frac12 \al_0, \beta - \frac12 \al_0) = 
(\beta, \beta) - (\beta, \al_0) + \frac14 (\al_0, \al_0) = (\beta, \beta) - 
\frac14 (\al_0, \al_0)$ since $1 = \pair{\beta}{\al_0} = 2 (\beta, 
\al_0)/(\al_0, \al_0)$ by (\ref{eq:pair}). Thus, $(\la, \la) > 0$ is 
determined by $(\beta, \beta)$, and for the latter there are at most two 
possible values.

To show the third property, pick a long weight $\la_0 \in \Phi$, i.e. $\la_0 
= \beta_0 - \frac12 \al_0$ for some long root $\beta_0 \in \Delta$ with 
$\pair{\beta_0}{\al_0} = 1$. Since our hypothesis implies that $\Delta$ is 
not of type $\mbox{C}_k$, such a $\beta_0$ and hence such a $\la_0$ exists.

Let $\gamma \in \Delta$ with $\pair \gamma {\al_0} = 1$, and let $\mu := 
\gamma - \frac12 \al_0 \in \Phi$. Then $\gamma \neq - \beta_0$ so that $\pair 
\gamma {\beta_0} \in \{-1, 0, 1, 2\}$, and $\pair \gamma {\beta_0} = 2$ iff 
$\gamma = \beta_0$ iff $\mu = \la_0$.

If $\pair \gamma {\beta_0} = 1$ then $\beta_0 - \gamma \in \Delta$ with 
$\pair{\beta_0 - \gamma}{\al_0} = 0$, so that $\la_0 - \mu = \beta_0 - \gamma
\in \Delta_\h$.

If $\pair \gamma {\beta_0} \in \{0, -1\}$ then $\pair \gamma {\al_0 - 
\beta_0} = 1 - \pair \gamma {\beta_0} \in \{1, 2\}$, thus when replacing 
$\la_0$ by $-\la_0$ and hence $\beta_0$ by $\al_0 - \beta_0$, then we can 
reduce to the previous cases.

{}From this description, it also follows that $\Phi_i = \{ \mu \in \Phi \mid 
\pair \mu {\beta_0} = \frac i2\}$

To show the fourth part, let $\Delta_\h^i := \{ \gamma \in \Delta_\h \mid 
\pair{\gamma}{\beta_0} = i\}$. Since $\pm \beta_0 \not \in \Delta_\h$, it 
follows that $\Delta_\h = \Delta_\h^{-1} \cup \Delta_\h^{0} \cup 
\Delta_\h^{1}$, and we let $\h_{\pm1} := \bigoplus_{\gamma \in \Delta_\h^{\pm 
1}} \g_\gamma$ and $\h_0 := \t_\h \oplus \bigoplus_{\gamma \in \Delta_\h^{0}} 
\g_\gamma$. Since $\Phi_i = \{ \mu \in \Phi \mid \pair \mu {\beta_0} = \frac 
i2\}$, the claims follow.

Finally, for the last part, note that by (\ref{eq:Adams}),
\[
(v_+ \circ v_-) w_r = (v_+ \circ w_r) v_- + 2 \om(v_-, w_r) v_+ + \om(v_+, 
w_r) v_- - \om(v_+, v_-) w_r.
\]
Now if $r > 0$ then $v_+ \circ w_r \in \h_{\frac32 + r} = 0$ and $\om(v_+, 
w_r) = 0$. Also, $\om(v_-, w_r) = 0$ for $r = 1/2$ showing the 
claim in this case, whereas for $r = 3/2$, $w_r$ is a scalar multiple of 
$v_+$ so that $\om(v_-, w_r) v_+ = \om(v_-, v_+) w_r$ which implies the 
assertion in this case as well. The proof of the cases $r < 0$ follows 
analogously.

Note that then for $w_r \in V_r, w_s \in V_s$ we also have $[v_+ \circ v_-, 
w_r \circ w_s] = ((v_+ \circ v_-) w_r) \circ w_s + w_r \circ ((v_+ \circ v_-) 
w_s) = -2(r+s) \om(v_+, v_-) w_r \circ w_s$, and the last assertion follows.
\end{proof}

Let us now suppose that $\h \subsetneq \sp(V, \om)$ and $\dim V > 2$, so that 
we have the decompositions from the lemma. Let $\psi \in \cR_\h^{(1)}$ be a weight 
element of weight $\mu \in \Phi$. Choose a long weight $\la_0 \in \Phi$, $\la_0 \neq \pm \mu$ 
so that -- after replacing $\la_0$ by its negative if necessary -- we may assume that 
$\mu \in \Phi_1$. Whence, $\psi(V_\la) \in \g_{\la + \mu}$ implies that $\psi(V_r) 
\subset \h_{r + \frac12}$ and, in particular, $\psi(V_{\frac32}) = 0$.

Note that $\g_{-\la_0 + \mu} = V_{\mu} \circ V_{-\la_0}$; namely, 
$\pair{-\la_0}{\la_0 - \mu} < 0$ so that $\g_{\la_0 - \mu} V_{-\la_0} = 
V_{-\mu}$ as all weight spaces are one dimensional. Thus, $(\g_{\la_0 - \mu}, 
V_{\mu} \circ V_{-\la_0}) = \om(\g_{\la_0 - \mu} V_{-\la_0}, V_\mu) = 
\om(V_{-\mu}, V_\mu) \neq 0$ so that $0 \neq V_{\mu} \circ V_{-\la_0} \subset 
\g_{-\la_0 + \mu}$ and the latter is one dimensional.

Pick $0 \neq v_{-\la_0} \in V_{-\la_0}$. Since $\psi(v_{-\la_0}) \in 
\g_{-\la_0 + \mu}$, there is a $u \in V_\mu$ such that $\psi(v_{-\la_0}) = u 
\circ v_{-\la_0}$. Therefore, after replacing $\psi$ by $\psi - \psi_u$, we 
may assume that $\psi(v_{-\la_0}) = 0$ and hence $\psi(V_{\pm\frac32}) = 0$.

If we let $v_\pm \in V_{\pm \frac32}$ with $\om(v_+, v_-) \neq 0$ and $w_\pm 
\in V_{\pm \frac12}$ 
then by (\ref{eq:defR1}) we must have
\be \label{eq:prolong1}
0 = R_{\psi(w_\pm)}(v_+, v_-) = 2 \om(v_+, v_-) \psi(w_\pm) + v_+ \circ 
(\psi(w_\pm) v_-) - v_- \circ (\psi(w_\pm) v_+).
\ee

Now $\psi(w_+) \in \h_1$, hence $\psi(w_+) v_+ = 0$ and thus
$
{}v_+ \circ (\psi(w_+) v_-) = [\psi(w_+), v_+ \circ v_-] = 2 \om(v_+, v_-) 
\psi(w_+),
$
where the last identity follows from the lemma. Then (\ref{eq:prolong1}) 
implies that $\psi(w_+) = 0$.

On the other hand, $\psi(w_-) \in \h_0$ so that $\psi(w_-) v_\pm \in 
V_{\pm \frac32}$, hence (\ref{eq:prolong1}) implies that $\psi(w_-) = c\ v_+ 
\circ v_-$ for some constant $c$. But then, $\psi(w_-) v_\pm = \mp 3 c\ 
\om(v_+, v_-) v_\pm$ by the lemma, and substituting into (\ref{eq:prolong1}) 
yields $c = 0$, i.e. $\psi(w_-) = 0$, and hence, $\psi = 0$.

Let $W \subset \cR_\h^{(1)}$ be the $H$-invariant complement of $\{\psi_u \mid 
u \in V\} \subset \cR_\h^{(1)}$, and let $\Psi$ be the set of weights of $W$. Since 
$W \subset \cR_\h^{(1)} \subset V \ot K(\h) \cong V \ot \h$, it follows that $\Psi 
\subset \Phi + \Delta_0$. Also, by what we have shown above, we must have 
$\Psi \cap \Phi = \emptyset$. 

Let $\nu \in \Psi$, and write $\nu = \mu + \al$ with $\mu \in \Phi$ and $\al \in 
\Delta_0$. Since $\Psi \cap \Phi = \emptyset$, it follows that that $\al \neq 0$, i.e., 
$\al \in \Delta$. If $\pair \mu \al < 0$ , then $\nu = \mu + \al \in \Psi \cap \Phi$; if 
$\pair \nu \al > 0$ then $\mu = \nu - \al \in \Psi \cap \Phi$, both of which are 
impossible. Thus, $2 = \pair \al \al = \pair \nu \al - \pair \mu \al \leq 0$ which is a 
contradiction.

Thus, we must have $\Psi = \emptyset$ and hence $W=0$.
\end{proof}

\noindent
Finally, we prove the following result which we shall need later on.

\begin{lem} \label{lem:reps} Let $\h \subset \sp(V, \om)$ be a special 
symplectic subalgebra, $\dim V \geq 4$, and let $\phi: V \ra V$ be a 
linear map such that
\be \label{Schur}
\phi(x) \circ y = \phi(y)\circ x \mbox{\ \ for all $x,y \in V$.}
\ee
Then $\phi$ is a multiple of the identity.
\end{lem}

\begin{proof} By (\ref{eq:Adams}) we have
\[
(\phi(x) \circ y) z - (\phi(x) \circ z) y = 2 \om(y,z) \phi(x) + \om(\phi(x), 
z) y - \om(\phi(x), y) z.
\]
But (\ref{Schur}) now implies that the cyclic sum in $x, y, z$ of the left
hand side vanishes, hence so does the cyclic sum of the right hand side, i.e.
\be \label{Schur2}
\ba{l}2 \left( \om(x,y) \phi(z) + \om(y,z) \phi(x) + \om(z,x) \phi(y)\right)\\
\\ \ \ \ \ =
(\om(\phi(y), z) - \om(\phi(z), y)) x + (\om(\phi(z), x) - \om(\phi(x), z)) y 
+ (\om(\phi(x), y) - \om(\phi(y), x)) z.
\ea \ee
For each $x \in V$, we may choose vectors $y,z \in V$ with $\om(x,y) = 
\om(x,z) = 0$ and $\om(y,z) \neq 0$ since $\dim V \geq 4$. Then 
(\ref{Schur2}) implies that $\phi(x) \in span(x,y,z)$ so that $\om(\phi(x), 
x) = 0$. Polarization then implies that $\om(\phi(x), y) + \om(\phi(y), x) = 
0$ for all $x, y \in V$.

Next, we take the symplectic form of (\ref{Schur2}) with $x$, and together 
with the preceding identity this yields
\[
\om(x,y) \om(\phi(x), z) = \om(x,z) \om(\phi(x), y)\ \ \ \ \mbox{for all 
$x,y,z \in V$.}
\]
Thus, $\om(x,y) \phi(x) = \om(\phi(x), y) x$ for all $x, y \in V$, and since 
for $0 \neq x \in V$ we can pick $y \in V$ such that $\om(x,y) \neq 0$, this 
implies that $\phi(x)$ is a scalar multiple of $x$ for all $x \in V$, whence 
$\phi$ is a multiple of the identity.
\end{proof}

\begin{keydf} \label{df:spsymplconn}
Let $(M, \om)$ be a (real or complex) symplectic manifold of (real or complex)
dimension at least $4$, equipped with a symplectic 
connection $\nabla$, i.e. a torsion free connection for which $\om$ is 
parallel. We say that $\nabla$ is a {\em special symplectic connection 
associated to the (simple) Lie group $\G$} if there is a special symplectic 
subgroup $\eH \subset \Sp(V, \om)$ associated to $\G$ in the sense of 
Definition~\ref{df:associated} such that the curvature of $\nabla$ is contained 
in $\cR_\h$ (cf. (\ref{eq:defRh}) and (\ref{eq:defRhspace})).
\end{keydf}

Definition~\ref{df:spsymplconn} coincides with the definition of 
special symplectic connections from the introduction. Namely, note that by 
the Ambrose-Singer holonomy theorem, the (restricted) holonomy of a special 
symplectic connection is evidently contained in $\eH \subset \Sp(V, \om)$, so 
that we have an $\eH$-reduction $B \ra M$ of the frame bundle of $M$ which is 
compatible with the connection.

If $\eH \subset \Sp(V, \om)$ is one of the subgroups $(i)$ or $(ii)$, then 
either there are two complementary parallel Lagrangian foliations 
(case $(i)$), or the connection is the Levi-Civita connection of a 
pseudo-K\"ahler metric (case $(ii)$). In either case, the condition that the 
curvature lies in ${\cal R}_\h$ is equivalent to the vanishing of the {\em 
Bochner curvature}, and such connections have been called {\em 
Bochner-bi-Lagrangian} in the first and {\em Bochner-K\"ahler} in the second 
case. For a detailed study of these connections, see \cite{Bochner}.

If $\eH = \Sp(V, \om)$ as in $(iii)$, then the condition that the curvature 
lies in ${\cal R}_\h$ is equivalent to saying that the connection is a (real 
or holomorphic) {\em symplectic connection of Ricci type} in the sense of 
\cite{Michel}.

Finally, if $\eH \subset \Sp(V, \om)$ is one of the subgroups $(iv) - 
(xviii)$ in Table~1, then, by Theorem~\ref{thm:curvatures}, {\em any} torsion 
free connection on such an $\eH$-structure must be special. In fact, these 
subgroups $\eH$ are precisely the absolutely irreducible proper subgroups of 
the symplectic group which can occur as the holonomy of a torsion free 
connection (cf. \cite{MS}, \cite{Habil}, \cite{Advances}).

It shall be the aim of the following sections to study special 
symplectic connections using the general algebraic setup established 
here rather than dealing with each of the geometric structures 
separately.

\section{Special symplectic connections and contact manifolds} 
\label{sec:construction}

We shall now recall some well known facts about contact manifolds and 
their symplectic reductions.

\begin{df} A {\em contact structure} on a real (complex, respectively) 
manifold $\cC$ is a smooth (holomorphic, respectively) distribution ${\cal D} 
\subset T\cC$ of codimension one such that the Lie bracket induces a 
non-degenerate map
\[
{\cal D} \times {\cal D} \longrightarrow T\cC/{\cal D} =: L.
\]
\end{df}
The line bundle $L \ra \cC$ is called the {\em contact line 
bundle}, and its dual can be embedded as
\be \label{eq:defL}
L^* = \left\{ \la \in T^*\cC \mid \la({\cal D}) = 0\right\} \subset 
T^*\cC.
\ee

Notice that we can define the line bundles $L \ra \cC$ and $L^* \ra \cC$ for 
an arbitrary distribution ${\cal D} \subset T\cC$ of codimension one. It is 
well known that such a distribution ${\cal D}$ yields a contact structure iff 
the restriction of the canonical symplectic form $\Om$ on $T^*\cC$ to $L^* 
\backslash 0$ is non-degenerate, so that in this case $L^* \backslash 0$ is a 
symplectic manifold in a canonical way. 

We regard $p: L^* \backslash 0 \ra \cC$ as a principal $(\R \backslash 
0)$-bundle ($\C^*$-bundle, respectively). In the real case, we may assume 
that $L^* \backslash 0$ has two components each of which is a 
principal $\R^+$-bundle, since this can always be achieved when replacing 
$\cC$ by a double cover if necessary. Thus, we get the principal 
$\R^+$-bundle ($\C^*$-bundle, respectively)
\[
p: \hat \cC \longrightarrow \cC,
\]
where $\hat \cC \subset L^* \backslash 0$ is a connected component.
The vector field $E_0 \in {\frak X}(\hat \cC)$ which generates the principal 
action is called {\em Euler field}, so that the flow along $E_0$ is 
fiberwise scalar multiplication in $\hat \cC \subset L^* \subset T^*\cC$. 
Thus, the {\em Liouville form} on $T^*\cC$ is given as $\la := E_0 \hook 
\Om$, and hence ${\frak L}_{E_0}(\Om) = \Om$ and $\Om = d\la$.
This process can be reverted. Namely, we have the following

\begin{prop} \label{prop:revertcontact}
Let $p: \hat \cC \ra \cC$ be a principal $\R^+$-bundle ($\C^*$-bundle, 
respectively) with a symplectic form $\Om$ on $\hat \cC$ such that 
${\frak L}_{E_0}\Om = \Om$ where $E_0 \in {\frak X}(\hat \cC)$ generates 
the principal action. Then there is a unique contact structure ${\cal D}$ on
$\cC$ and an equivariant imbedding $\imath: \hat \cC \hookrightarrow L^* 
\backslash 0 \subset T^*\cC$ with $L^*$ from (\ref{eq:defL}) such that $\Om$ 
is the pullback of the canonical symplectic form on $T^*\cC$ to $\hat \cC$.
\end{prop}

\begin{proof}
By hypothesis, $\Om = d\la$ where $\la := (E_0 \hook \Om)$. Since $\la(E_0) =
0$, there is for each $x \in \hat \cC$ a unique ${\und \la}_x \in T^*_{p(x)} 
\cC$ satisfying $p^*({\und \la}_x) = \la_x$. Moreover, ${\frak L}_{E_0}(\la) 
= \la$, hence ${\und \la}_{e^t x} = e^t {\und \la}_x$ for all $t 
\in \F$, so that the codimension one distribution ${\cal D} := dp(\ker (\la)) 
\subset T\cC$ is well defined, and the correspondence $x \mapsto \und \la_x$ 
yields an equivariant imbedding $\hat \cC \hookrightarrow L^* \backslash 0$ 
whose image is thus a connected component of $L^* \backslash 0$.
Moreover, by construction, $\la$ is the restriction of the Liouville form to 
$\hat \cC \subset L^* \backslash 0 \subset T^* \cC$. Since $\Om = d\la$ is 
non-degenerate on $\hat \cC$ by assumption, it follows that ${\cal D}$ is a 
contact structure.  
\end{proof}

\noindent
Next, we define the fiber bundle 
\[
{\frak R} := \{ (\la, \hat \xi) \in \hat \cC \times T\hat \cC \subset T^*\cC 
\times T\hat \cC \mid \la(dp(\hat \xi)) = 1\}.
\]
Projection onto the first factor yields a fibration ${\frak R} \ra 
\hat \cC$ whose fiber is an affine space.

We call a vector field $\xi$ on $\cC$ a {\em contact symmetry} 
if ${\frak L}_\xi (\cal D) \subset {\cal D}$. This means that the flow 
along $\xi$ preserves the contact structure ${\cal D}$. For each 
contact symmetry $\xi$ on $\cC$, there is a unique vector field $\hat 
\xi \in {\frak X}(\hat \cC)$, called the {\em Hamiltonian lift of $\xi$},
satisfying $dp(\hat \xi) = \xi$ and ${\frak L}_{\hat \xi}\la = 0$, so that 
${\frak L}_{\hat \xi}\Om = 0$.

We call $\xi$ a {\em transversal} contact symmetry if in addition $\xi 
\not\in {\cal D}$ at all points. Equivalently, we have $\Om(E_0, \hat 
\xi) \neq 0$ everywhere. In the real case, we say that $\xi$ is {\em 
positively transversal} if $\Om(E_0, \hat \xi) > 0$ everywhere, while in the
complex case it is convenient to call {\em any} transversal vector field 
positively transversal.

Given a positively transversal contact symmetry $\xi$ with Hamiltonian lift
$\hat \xi$, there is a unique 
section $\la$ of the bundle $p: \hat \cC \ra \cC$ such that $\la(\xi) \equiv 
1$, and hence we obtain a section of the bundle ${\frak R} \ra \hat \cC \ra 
\cC$
\be \label{eq:sectionR}
\sigma_\xi: \cC \longrightarrow {\frak R},\ \ \ \ \ \sigma_\xi := (\la, \hat
\xi) \in {\frak R}.
\ee

We call an open subset $U \subset \cC$ {\em regular} w.r.t. the 
transversal contact symmetry $\xi$ if there is a submersion $\pi_U: U \ra 
M_U$ onto some manifold $M_U$ whose fibers are connected lines tangent 
to $\xi$. Evidently, since $\xi$ is pointwise non-vanishing, $\cC$ 
can be covered by regular open subsets. 

Since $\xi$ is a contact symmetry, it follows that $\xi \hook d\la = 0$ and 
${\frak L}_\xi \la = 0$. Thus, on each $M_U$ there is a unique symplectic form
$\om$ such that 
\be \label{eq:symplMU}
\pi_U^*\om = -2 d\la,
\ee
where the factor $-2$ only occurs to make this form coincide with one we shall
construct later on.

To link all of this to our situation, let $\g$ be a $2$-gradable simple real 
or complex Lie algebra and let $\G$ be the corresponding connected Lie group 
with trivial center $Z(\G) = \{1\}$. Recall the decomposition
\[
\g = \g^{-2} \oplus \g^{-1} \oplus \g^0 \oplus \g^1 \oplus \g^2 
\cong \F e_-^2 \oplus (e_- \ot V) \oplus (\F e_+ e_- \oplus \h) \oplus 
(e_+ \ot V) \oplus \F e_+^2
\]
from (\ref{eq:decomposeg}). We let $\mu := g^{-1}dg$ be the left invariant 
Maurer-Cartan form on $\G$, which we can decompose as 
\be \label{eq:decomposemu}
\mu = \sum_{i=-2}^2 \mu_i,\ \ \ \ \ \mu_0 = \mu_\h + \nu_0 e_+ e_-
\ee 
where $\mu_i \in \Om^1(\G) \ot \g^i$, $\mu_\h \in \Om^1(\G) \ot \h$ 
and $\nu_0 \in \Om^1(\G)$.
Furthermore, we define the subalgebras
\[ 
\ba{lll}
\p := \g^0 \oplus \g^1 \oplus \g^2, & \mbox{\ \ \ and\ \ \ } & 
\p_0 := \h \oplus \g^1 \oplus \g^2,
\ea
\]
and we let $\eP, \eP_0 \subset \G$ be the corresponding connected subgroups. 
Using the bilinear form $(\ ,\ )$ from (\ref{eq:Killing}), we identify $\g$ 
and $\g^*$, and recall the root cone from (\ref{eq:cones}) and its 
(oriented) projectivization
\be \label{eq:projectivize}
\hat \cC := \G \cdot e_+^2 \subset \g \cong \g^*,\ \ \ \ \ \ 
\cC := p(\hat \cC) \subset \P^o(\g) \cong \P^o(\g^*),
\ee
where $\P^o(\g)$ is the set of {\em oriented} lines in $\g$, i.e. $\P^o 
\cong S^d$ if $\F = \R$, and $\P^o \cong \C\P^d$ if $\F = \C$, where $d = 
\dim \g - 1$, and where $p: \g \backslash 0 \ra \P^o(\g)$ is the principal 
$\R^+$-bundle ($\C^*$-bundle, respectively) defined by the canonical 
projection. Thus, the restriction $p: \hat \cC \ra \cC$ is a principal
bundle as well.

Being a coadjoint orbit, $\hat \cC$ carries a canonical $\G$-invariant 
symplectic structure $\Om$. Moreover, the {\em Euler vector field} 
defined by
\[
E_0 \in {\frak X}(\hat \cC),\ \ \ \ \ (E_0)_{v} := v
\]
generates the principal action of $p$ and satisfies ${\frak L}_{E_0}(\Om) = 
\Om$, so that the distribution ${\cal D} = dp(E_0^{\perp_\Om}) \subset T\cC$ 
yields a $\G$-invariant contact distribution on $\cC$ by 
Proposition~\ref{prop:revertcontact}.

\begin{lem} As homogeneous spaces, we have $\cC = \G/\eP$, $\hat \cC = 
\G/\eP_0$ and ${\frak R} = \G/\eH$. Moreover, the fiber bundles ${\frak R} 
\ra \hat \cC \ra \cC$ from before are equivalent to the corresponding 
homogeneous fibrations.
\end{lem}

\begin{proof} The first two statements follow immediately from (\ref{eq:projectivize}). 
Using the pairing $(\ ,\ )$ to identify $\g$ and $\g^*$, it follows that the fiber of 
${\frak R}$ over $e_+^2 \in \hat \cC$ can be identified with
\[
{\frak R}_{e_+^2} = \left\{\frac12 e_-^2 + e_- \ot v + t e_+ e_- + \p_0 \mid v \in V, t \in 
\F\right\} \subset \g/\p_0 \cong T_{e_+^2} \hat \cC.
\]
Now it is straightforward to verify that $\eP_0 = \exp(\p_0)$ acts transitively 
on this set. Moreover, for all $p_0 \in \p_0$ one calculates that 
$(\ad(\frac12 e_-^2 + p_0))^2 (e_+^2) \in \F (\frac12 e_-^2 + p_0)$ iff $p_0 
= 0$. Since $(\ad_x)^2(\g) \subset \F x$ for all $x \in \hat \cC$, it follows 
that $(\frac12 e_-^2 + \p_0) \cap \hat \cC = \frac12 e_-^2$, and hence each 
of the cosets $\{\frac12 e_-^2 + e_- \ot v + t e_+ e_- + \p_0\} \in \g/\p_0$ 
has a unique representative in $\hat \cC$. 

{}From all of this it now follows that $\G$ acts transitively on ${\frak 
R}$, and the stabilizer of the pair $(e_+^2, \frac12 e_-^2 + \p_0)$ equals 
the stabilizer of the pair $(e_+^2, \frac 12 e_-^2)$ which is $\eH$ by 
Proposition~\ref{prop:Hconnected} as $Z(\G) = \{1\}$. Thus, ${\frak R} = 
\G/\eH$ as claimed.

The fibers of the homogeneous fibrations ${\frak R} \ra \hat \cC$ and ${\frak 
R} \ra \cC$ are connected, and since ${\frak R} = \G/\eH$ and $\eH$ is 
connected, it follows that the stabilizers of $e_+^2 \in \hat \cC$ and 
$[e_+^2] \in \cC$ are connected as well. Since the Lie algebras of these 
stabilizers are evidently $\p_0$ and $\p$, respectively, the claim follows.
\end{proof}

For each $a \in \g$ we define the vector fields $a^* \in {\frak 
X}(\cC)$ and $\hat a^* \in {\frak X}(\hat \cC)$
corresponding to the infinitesimal action of $a$, i.e.
\be \label{eq:xi0}
(a^*)_{[v]} := \left. \frac d{dt} \right|_{t=0} (\exp(t a) \cdot [v])
\ \ \ \ \ \mbox{and}\ \ \ \ 
(\hat a^*)_v := \left. \frac d{dt} \right|_{t=0} (\exp(t a) \cdot v).
\ee
Note that $a^*$ is a contact symmetry and $\hat a^*$ is its 
Hamiltonian lift. Let
\be \label{eq:cC0}
\hat \cC_a := \{ \la \in \hat \cC \mid \la(a^*) \in \R^+ (\in 
\C^*, \mbox{ respectively})\}\ \ \ \ \ 
\mbox{and}\ \ \ \ \ \cC_a := p(\hat \cC_a) \subset \cC,
\ee
so that $p: \hat \cC_a \ra \cC_a$ is a principal $\R^+$-bundle 
($\C^*$-bundle, respectively) and the restriction of $a^*$ to $\cC_a$ is 
a positively transversal contact symmetry. Therefore, we obtain the section 
$\sigma_a: \cC_a \ra {\frak R} = \G/\eH$ from (\ref{eq:sectionR}).

Let $\pi: \G \ra \G/\eH = {\frak R}$ be the canonical projection, and let 
$\Gamma_a:= \pi^{-1}(\sigma_a(\cC_a)) \subset \G$. Then evidently, the 
restriction $\pi: \Gamma_a \ra \sigma_a(\cC_a) \cong \cC_a$ is a 
(right) principal $\eH$-bundle.

\begin{thm} \label{thm:Maurer-Cartan}
Let $a \in \g$ be such that $\cC_a \subset \cC$ from (\ref{eq:cC0}) is 
non-empty, define $a^* \in {\frak X}(\cC)$ and $\hat a^* \in 
{\frak X}(\hat \cC)$ as in (\ref{eq:xi0}), and let $\pi: \Gamma_a \ra \cC_a$ 
with $\Gamma_a \subset \G$ be the principal $\eH$-bundle from above. Then 
there are functions $\rho : \Gamma_a \ra \h$, $u : \Gamma_a \ra V$, $f: 
\Gamma_a \ra \F$ such that 
\be \label{eq:formA}
\Ad_{g^{-1}}(a) = \frac12 e_-^2 + \rho + e_+ \ot u + \frac12 f e_+^2
\ee 
for all $g \in \Gamma_a$. Moreover, the restriction of the components 
$\mu_{-2} + \mu_{-1} + \mu_\h$ of the Maurer-Cartan form 
(\ref{eq:decomposemu}) to $\Gamma_a$ yields a pointwise linear isomorphism 
$T\Gamma_a \ra \h \oplus \g^{-1} \oplus \g^{-2}$, and if we decompose this 
coframe as
\[
\mu_{-2} + \mu_{-1} + \mu_\h = -2 \kappa\ \left(\frac12 e_-^2 + 
\rho\right) + e_- \ot \th + \eta,\ \ \ \ 
\kappa \in \Om^1(\Gamma_a),\ \ \th \in \Om^1(\Gamma_a) \ot V,\ \ \ 
\eta \in \Om^1(\Gamma_a) \ot \h,
\]
then $\kappa = -\frac12 \pi^*(\la)$ where $\la \in \Om^1(\cC_a)$ is the 
contact form for which $\sigma_a = (\la, \hat a^*)$. Moreover, we have the 
structure equations
\be \label{eq:structurekappa}
d\kappa = \frac12 \om(\th \w \th),
\ee
and
\be \label{eq:structureGamma} \ba{lll}
\ba{lll}
d\th + \eta \w \th & = & 0,\\ \\
d\eta + \frac12 [\eta, \eta] & = & R_{\rho} (\th \w \th),
\ea
& \mbox{\hspace{2cm}} & 
\ba{lll}
d\rho + [\eta, \rho] & = & u \circ \th\\ \\ 
du + \eta \cdot u & = & (\rho^2 + f) \cdot \th\\ \\ 
df + d(\rho, \rho) & = & 0.
\ea \ea \ee
\end{thm}

\begin{proof}
According to the above identifications, we have $g \in \Gamma_a$ iff $(g 
\cdot e_+^2, g \cdot (\frac 12 e_-^2 + \p_0)) = \sigma_a([g \cdot e_+^2])$ 
iff $g \cdot (\frac 12 e_-^2 + \p_0) = (\hat a^*)_{g \cdot e_+^2}$ iff 
$(\Ad_{g^{-1}}(\hat a^*))_{e_+^2} = \frac 12 e_-^2 \mod \p_0$ iff 
$\Ad_{g^{-1}}(a) = \frac 12 e_-^2 \mod \p_0$, i.e.
\be \label{eq:Gamma0}
\ba{ll} & 
\Gamma_a = \left\{ g \in \G \ \left|\ \Ad_{g^{-1}}(a) \in Q\right. 
\right\},\\ \mbox{where}\\ 
& Q := \frac12 e_-^2 + \p_0 = \left\{ \frac12 e_-^2 + \rho + e_+ \ot u + 
\frac12 f e_+^2 \mid \rho \in \h, u \in V, f \in \F \right\},
\ea
\ee
and from this (\ref {eq:formA}) follows.
Thus, if $dL_g v \in T_g \Gamma_a$ with $v \in \g$, then we must have 
\[
\p_0 \ni \left. \frac d{dt} \right|_{t=0} \left(\Ad_{(g \exp(t 
v))^{-1}}(a)\right) = -[v, \Ad_{g^{-1}}(a)] = - \left[v, \frac12 e_-^2 + 
\rho + e_+ \ot u + \frac12 f e_+^2\right],
\]
and from here it follows by a straightforward calculation that $v$ 
must be contained in the space
\be \label{eq:TGamma}
\F \Ad_{g^{-1}} a \oplus \left\{\left. \ e_- \ot x + e_+ \ot \rho x + 
\frac12 \om(u,x) e_+^2\ \right|\  x \in V \right\} \oplus \h,
\ee
and since $v$ was arbitrary, it follows that $\mu(T_g \Gamma_a)$ is contained 
in (\ref{eq:TGamma}). In fact, a dimension count yields that 
$\dim (\mu(T_g \Gamma_a)) = \dim \Gamma_a = \dim \cC_a + \dim \eH$ coincides 
with the dimension of (\ref{eq:TGamma}), hence (\ref{eq:TGamma}) 
equals $\mu(T_g \Gamma)$, i.e. $\mu_{-2} + \mu_{-1} + \mu_\h: T\Gamma_a \ra 
\g^{-2} \oplus \g^{-1} \oplus \h$ yields a pointwise isomorphism. 

Let $\xi_a$ denote the {\em right} invariant vector field on $\G$ characterized by 
$\mu(\xi_a) = \Ad_{g^{-1}}(a)$. Then the flow of $\xi_a$ is {\em left} multiplication 
by $\exp(t a)$ and hence evidently leaves $\Gamma_a$ invariant. Moreover, by 
(\ref{eq:TGamma}) we have
\be \label{eq:properties-xi_a} \ba{lll}
{\frak L}_{\xi_a}^*(\mu) = 0, & \mbox{ and } & d\rho(\xi_a) = du(\xi_a) = df(\xi_a) = 0.
\ea \ee
Let us write the components of the Maurer-Cartan form $\mu$ as
\[ \ba{lll}
\mu_{\pm2} := \kappa_\pm e_\pm^2, & \mu_{\pm1} := e_\pm \ot \al_\pm, & \mu_0 := 
\nu e_+ e_- + \mu_\h
\ea\]
with $\kappa_\pm, \nu \in \Om^1(G)$, $\al_\pm \in \Om^1(G) \ot V$ and $\mu_\h 
\in \Om^1(G) \ot \h$. Now $\mu(T\Gamma_a)$ is given by (\ref{eq:TGamma}) so that 
by  (\ref{eq:formA}), the restriction of $\mu$ to $\Gamma_a$ satisfies
\[ \ba{llll}
\nu = 0, & \al_+ = \rho \al_- - 2 u \kappa, & \kappa_+ = \frac12 \om(u, \th) - 2 f \kappa & 
\mu_\h = \eta - 2 \kappa \rho,
\ea \]
where $\kappa := - \kappa_-$, $\th := \al_-$ and $\eta := \mu_\h + 2 \kappa \rho$. 
Substituting this into the Maurer-Cartan equation $d\mu + \frac12 [\mu, \mu] = 0$, a 
straightforward calculations yields (\ref{eq:structurekappa}) and 
\be \label{eq:MC-2} \ba{rll}
d\th + \eta \w \th & = & 0,\\ \\
d\eta + \frac12 [\eta, \eta] - R_{\rho} (\th \w \th) & = & -2 \kappa \w (d\rho + [\eta, \rho] - 
u \circ \th),\\Ê\\
(d\rho + [\eta, \rho] - u \circ \th) \w \th & = & -2 \kappa \w (du + \eta \cdot u - (\rho^2 + f) 
\cdot \th),\\ \\ 
\om(du + \eta \cdot u - (\rho^2 + f) \cdot \th, \th) & = & -2 \kappa \w (df + d(\rho, \rho)).
\ea \ee

By (\ref{eq:properties-xi_a}), we have $\th(\xi_a) = \eta(\xi_a) = 0$, $\kappa(\xi_a) \equiv 
-\frac12$, and $\xi_a \hook d\th = \xi_a \hook d\eta = 0$. Thus, the contraction of the left 
hand sides of (\ref{eq:MC-2}) with $\xi_a$ vanishes, and from there, 
(\ref{eq:structureGamma}) follows.

Note that $dp(\xi_a) = \hat a^*$, where $p: \Gamma_a \ra \hat \cC$ is the canonical 
projection, and from (\ref {eq:formA}) it follows that $\la(a^*) = -2 \kappa(\tilde a^*) \equiv 
1$, so that $(\la, \hat a^*) \in {\frak R}$ which shows the final assertion.
\end{proof}

With these structure equations, we are now ready to prove the following 
result which immediately implies Theorem~A of the introduction.

\begin{thm} \label{thm:canonicalconn}
Let $a \in \g$ and $\cC_a \subset \cC$ as before. Let $U \subset \cC_a$ be a
regular open subset , i.e. the local quotient $M_U := \T_a^{loc} \backslash 
U$ is a manifold, where
\[
\T_a := \exp(\F a) \subset \G.
\]
Let $\om \in \Om^2(M)$ be the symplectic form from (\ref{eq:symplMU}). Then 
$M_U$ carries a canonical special symplectic connection associated to $\g$,
and the (local) principal $\T_a$-bundle $\pi: U \ra M$ admits a connection 
$\kappa \in \Om^1(U)$ whose curvature is given by $d\kappa = \pi^*(\om)$.
\end{thm}

\begin{proof}
Let us consider the commutative diagram
\be \label{diagramGamma}
\xymatrix{
\Gamma_a \ar[r]^{\mbox{\footnotesize $\T_a$}} 
\ar[d]_{\mbox{\footnotesize $\eH$}} & \T_a \backslash \Gamma_a 
\ar[d]_{\mbox{\footnotesize $\eH$}}\\
\cC_a \ar[r]^{\mbox{\footnotesize $\T_a$}} & \T_a \backslash \cC_a}
\ee
where the maps $\pi: \Gamma_a \ra \T_a \backslash \Gamma_a$ and $\Gamma_a \ra 
\cC_a$ are principal bundles with the indicated structure groups, 
whereas the arrows $\T_a \backslash \Gamma_a \ra \T_a \backslash \cC_a$ and 
$\cC_a \ra \T_a \backslash \cC_a$ stand for fibrations with a locally free, 
but not necessarily free group action of the indicated structure group.

It follows now immediately from (\ref{eq:structurekappa}) and 
(\ref{eq:structureGamma}) that $\th + \eta$ and $\kappa$ are the pull backs of
one forms on $\T_a \backslash \Gamma_a$ and $\cC_a$, respectively, and we 
shall by abuse of notation denote these forms by the same symbols.

Let $U \subset \cC_a$ be a regular open subset, let $\Gamma_U:= \pi^{-1}(U)
\subset \Gamma_a$ and $B := \T_a^{loc} \backslash \Gamma_U$ be the 
corresponding subsets. It follows then that the induced commutative diagram

\[
\xymatrix{
\Gamma_U \ar[r]^{\mbox{\footnotesize $\T_a^{loc}$}} 
\ar[d]_{\mbox{\footnotesize $\eH$}} & B 
\ar[d]_{\mbox{\footnotesize $\eH$}}\\
U \ar[r]^{\mbox{\footnotesize $\T_a^{loc}$}} & M
}
\]
consists of {\em (local) principal} bundles, and $B$ and $U$ carry a $V 
\oplus \h$-valued coframe $\th + \eta$ and a one form $\kappa$, respectively,
satisfying $d\kappa = \pi^*(\om)$ and (\ref{eq:structureGamma}), where 
$\om \in \Om^2(M)$ is the canonically induced symplectic form from 
(\ref{eq:symplMU}).

Standard arguments now show that $B \ra M$ is an $\eH$-structure with 
tautological one form $\th$, and $\eta$ defines a connection on $M$. By 
(\ref{eq:structureGamma}), this connection is torsion free and its curvature 
is given by $R_{\rho} (\th \w \th)$, i.e. this connection is special 
symplectic in the sense of Definition~\ref{df:spsymplconn}. 
\end{proof}

\begin{rem} \label{rem:conjugate}
{\em If we replace $a$ by $a' := \Ad_{g_0}(a)$, then it is clear that in the 
above construction we have $\Gamma_{a'} = L_{g_0} \Gamma_a$. Thus, 
identifying $\Gamma_a$ and $\Gamma_{a'}$ via $L_{g_0}$, the functions $\rho + 
\mu + f$ and the forms $\kappa + \th + \om$ will be canonically identified 
and hence both satisfy (\ref{eq:structureGamma}). Therefore, the connections 
from the preceding theorem only depend on the adjoint orbit of $a$.

Also, let $e^{t_0}$ with $t_0 \in \F$. Since $\cC_a = \cC_{e^{t_0} a}$ and 
$\T_a = \T_{e^{t_0} a}$, the above construction yields equivalent connections 
when replacing $a$ by $e^{t_0} a$. In this case, however, the symplectic form 
$\om$ on the quotient will be replaced by $e^{-t_0} \om$.
}
\end{rem}

\section{The structure equations}\label{sec:develope}

In this section, we shall revert the process of the preceding section, 
showing that any special symplectic connection is equivalent to the 
ones given in Theorem~\ref{thm:canonicalconn} in a sense which is to 
be made precise. We begin by deriving the structure equations for
special symplectic connections.

\begin{prop} \label{prop:allareinduced}
Let $(M, \om, \nabla)$ be a (real or complex) symplectic manifold of dimension
$\geq 4$ with a special symplectic connection of regularity $C^4$ associated 
to the Lie algebra $\g$, and let $\h \subset \g$ be as before. Then there is 
an associated $\tilde \eH$-structure $\pi: B \ra M$ on $M$ which is 
compatible with $\nabla$, where $\tilde \eH \subset \Sp(V, \om)$ is a Lie 
subgroup with Lie algebra $\h$, and there are maps $\rho: B \ra \h$, $u: B 
\ra V$ and $f: B \ra \F$, where $\F = \R$ or $\C$, such that the tautological 
form $\th \in \Om^1(B) \ot V$, the connection form $\eta \in \Om^1(B) \ot \h$ 
and the functions $\rho, u$ and $f$ satisfy the structure equations 
(\ref{eq:structureGamma}).
\end{prop}

To slightly simplify our arguments, we shall assume that $\tilde \eH = \eH$ 
is connected, which can be achieved by passing to an appropriate covering of
$M$. However, our results (and in particular Theorem~B) also hold if $\tilde 
\eH$ is {\em not} connected.

For clarification, we restate the structure equations 
(\ref{eq:structureGamma}) as follows. If for $h \in \h$ and $x \in V$ we let 
the vector fields $\xi_h, \xi_x \in \frak{X}(B)$ be the vector fields which 
are characterized by
\be \label{eq:tautologicalfields}
\th(\xi_h) \equiv 0,\ \ \eta(\xi_h) \equiv h\ \ \ \ 
\mbox{and}\ \ \ \ \th(\xi_x) \equiv x,\ \ \eta(\xi_x) \equiv 0,
\ee
then for all $h, l \in \h$ and $x,y \in V$,
\be \label{eq:structure2}
\ba{llllllll}
{}[\xi_h, \xi_l] & = \xi_{[h,l]}, & \ \ & [\xi_h, \xi_x] & = 
\xi_{hx}, & \ \ & [\xi_x, \xi_y] & = - 2 \om(x,y) \xi_\rho - \xi_{x 
\circ \rho y} + \xi_{y \circ \rho x}\\ \\
\xi_h(\rho) & = -[h, \rho], & & \xi_h(u) & = -h u, & & \xi_h(f) & = 0,\\ \\
\xi_x(\rho) & = u \circ x, & & \xi_x(u) & = (\rho^2 + f) x, & & \xi_x(f) 
& = -2 \om(\rho u, x)
\ea
\ee

The proof can be found e.g. in \cite{Michel} for the case of connections of 
Ricci type, in \cite{Advances} for the case of the special symplectic 
holonomies and in \cite{Bochner} in the case of Bochner K\"ahler 
metrics. But for the sake of completeness (and since our notation here is 
slightly different) we restate it here.

\begin{proof} Let $F$ be the $\eH$-structure on the manifold $M$, and denote 
the tautological and the connection $1$-form on $F$ by $\th$ and $\eta$, 
respectively. Since by hypothesis, the curvature maps are all contained in 
$\cR_\h$, it follows that there is an $\eH$-equivariant map $\rho: B \ra \h$ 
such that the curvature at each point is given by $R_\rho$ with the notation 
from (\ref{eq:defRh}). Thus, we have the structure equations
\be \label{eq:structuresimple} 
\ba{ll}
d\th + \eta \w \th & = 0\\
d\eta + \frac12 [\eta, \eta] & = R_\rho \cdot (\th \w \th),
\ea \ee
The $\eH$-equivariance of $\rho$ yields that $\xi_h(\rho) = -[h, \rho]$ for 
all $h \in \h$. Moreover, since the covariant derivative of the curvature is 
represented by $\xi_x(\rho)$ for all $x \in V$ and this must lie in 
$\cR^{(1)}_\h$, it follows by Proposition~\ref{prop:R1h} that $\xi_x(\rho) = 
u \circ \rho$ for some $\eH$-equivariant map $u: B \ra V$, which shows the 
asserted formula 
\be \label{eq:struct4}
d\rho + [\eta, \rho] = u \circ \th.
\ee

Since $u$ is $\eH$-equivariant, it follows that $\xi_h(u) = -hu$ for all $h 
\in \h$. Also, differentiation of (\ref{eq:struct4}) yields that for all $x, 
y \in V$
\[
\left(\xi_x u - \rho^2 x\right) \circ y = \left(\xi_y u - \rho^2 y\right) 
\circ x.
\]
Thus, by Lemma~\ref{lem:reps} it follows that there is a smooth function $f: B
\ra \F$ for which $\xi_x u - \rho^2 x = f x$ for all $x \in V$ so that
\be \label{eq:struct5}
du + \eta \cdot u = (\rho^2 + f) \th.
\ee
Finally, taking the exterior derivative of (\ref{eq:struct5}) yields that $df 
+ d(\rho, \rho) = 0$.
\end{proof}

It is now our aim to construct the equivalent to the line bundle 
$\Gamma \ra B$ from the preceding section. Motivated by (\ref{eq:Gamma0}) and
(\ref{eq:TGamma}), we define the following function $A$ and one form $\sigma$

\be \label{eq:defsigma}
\ba{lll}
A: B \longrightarrow Q \subset \g, & \mbox{\ \ \ \ } & 
A := \frac12 e_-^2 + \rho + e_+ \ot u + \frac12 f e_+^2,\\ \\
\sigma \in \Om^1(B) \ot \g, & & \sigma := e_- \ot \th + \eta + e_+ 
\ot (\rho \th) + \frac12 \om(u, \th) e_+^2,
\ea \ee
where $Q := \frac12 e_-^2 + \p_0 \subset \g$ is the affine hyperplane from 
(\ref{eq:Gamma0}). It is then straightforward to verify that 
(\ref{eq:structureGamma}) is equivalent to 

\be \label{eq:dA}
dA = -[\sigma, A]\ \ \ \ \ \ \mbox{and}\ \ \ \ \ \ d\sigma + 
\frac12 [\sigma, \sigma] = 2 \pi^*(\om) A.
\ee
Let us now enlarge the principal $\eH$-bundle $B \ra M$ to the principal
$\G$-bundle
\[
\B := B \times_\eH \G \longrightarrow M,
\]
where $\eH$ acts on $B \times \G$ from the right by $(b, g) \cdot h := (b 
\cdot h, h^{-1} g)$, using the principal $\eH$-action on $B$ in the first
component. Evidently, the inclusion $B \times \eH \hookrightarrow B \times \G$
induces an embedding $B \hookrightarrow \B$.

\begin{prop} \label{prop:extendbundle}
The function $\A$ and the one form $\al$ defined by
\be \label{eq:def-A-alpha}
\ba{lllll}
\A: \B \longrightarrow \g, & \mbox{\ \ \ \ } & \A[(b,g)] := 
\Ad_{g^{-1}}(A(b)),\\ \\ 
\al \in \Om^1(\B) \ot \g, & & \al_{[(b,g)]} := \Ad_{g^{-1}} \sigma_b + \mu,
\ea \ee
on $\B$ are well defined, where $\mu = g^{-1} dg \in \Om^1(\G) \ot \g$ is the 
left invariant Maurer-Cartan form on $\G$, and the restriction of $\A$ to $B 
\subset \B$ coincides with $A$. Moreover, $\al$ yields a connection on the 
principal $\G$-bundle $\B \ra M$ which satisfies
\be \label{eq:dA-hat}
d\A = -[\al, \A]\ \ \ \ \ \mbox{and}\ \ \ \ \ d\al + \frac12 [\al, \al] = 2 
\pi^*(\om) \A.
\ee
\end{prop}

\begin{proof}
First, note that $A: B \ra \eH$ and $\sigma \in \Om^1(B) \ot \g$ are 
$\eH$-equivariant, i.e. $R_h^* A = \Ad_{h^{-1}} A$ and $R_h^* \sigma = 
\Ad_{h^{-1}} \sigma$. Thus, if we define the function $\hat \A$ and the one form
$\hat \al$ by
\[ \ba{ll}
\hat \A := \Ad_{g^{-1}}(A): & B \times \G \longrightarrow \g\\ \\
\hat \al := \Ad_{g^{-1}} \sigma + \mu & \in \Om^1(B \times \G) \ot \g, 
\ea \]
then $\hat \A(bh, h^{-1}g) = \hat \A(b,g)$, so that $\hat \A$ is the pull 
back of a well defined function $\A: \B \ra \g$. Also, $\hat \al$ is 
invariant under the right $\eH$-action from above, and for $h \in \h$ we have
\[
\hat \al((\xi_h)_b, dR_g(-h)) =  \Ad_{g^{-1}} (\sigma_b(\xi_h)) - 
\mu(dR_g(h)) = \Ad_{g^{-1}} (h) - \Ad_{g^{-1}}(h) = 0,
\]
so that $\hat \al$ is indeed the pull back of a well defined form $\al \in 
\Om^1(\B) \ot \g$. Moreover, $R_g^*(\hat \al) = \Ad_g^{-1} \hat \al$ is easily
verified, and since $\hat \al$ coincides with $\mu$ on the fibers of the 
projection $B \times \G \ra B$, it follows that the value of $\hat \al$ on 
each left invariant vector field on $\G$ is constant. Since the left 
invariant vector fields generate the principal right action of the bundle $B 
\times \G \ra B$, it follows that $\hat \al$ is a connection on this bundle, 
hence so is $\al$ on the quotient $\B \ra M$.

Finally, to show (\ref{eq:dA-hat}) it suffices to show the corresponding 
equations for $\hat \al$ and $\hat \A$. We have
\[
d\hat \A = -[\mu, \Ad_{g^{-1}}(A)] + \Ad_{g^{-1}}(dA) = -[\mu, \hat \A] -
\Ad_{g^{-1}}([\sigma, A]) = -[\mu, \hat \A] - [\Ad_{g^{-1}} \sigma, \hat \A] =
-[\hat \al, \hat \A]
\]
by (\ref{eq:dA}), and
\[
\ba{lll}
d\hat \al + \frac12 [\hat \al, \hat \al] & = & (-[\mu, \Ad_{g^{-1}} \sigma] + 
\Ad_{g^{-1}} d\sigma + d\mu) + \frac12 (\Ad_{g^{-1}} [\sigma, \sigma] + 2 
[\mu, \Ad_{g^{-1}} \sigma] + [\mu, \mu])\\ \\
& = & \Ad_{g^{-1}} (d\sigma + \frac12 [\sigma, \sigma]) + d\mu + \frac12 [\mu,
\mu]\\ \\
& = & \Ad_{g^{-1}} (2 \pi^*(\om) A) = 2 \pi^*(\om) \hat \A,
\ea
\]
where the second to last equation follows from the Maurer-Cartan equation and 
(\ref{eq:dA}).
\end{proof}

\noindent
{\bf Proof of Theorem~B.} 
Let $\hat M \subset \B$ be a holonomy reduction of $\al$, and let $\hat \T
\subset \G$ be the holonomy group, so that the restriction $\hat M \ra M$
becomes a principal $\hat \T$-bundle. By the first equation
of (\ref{eq:dA-hat}), it follows that $\hat M \subset \A^{-1}(a)$ for some
$a \in \g$, and by choosing the holonomy reduction such that it contains an
element of $B \subset \B$, we may assume w.l.o.g. that $a \in Q$. We let
\be \label{eq:defS}
\hat \S := Stab(a) = \{ g \in \G \mid \Ad_g a = a\} \subset \G\ \ \ \ 
\mbox{and}\ \ \ \ \hat \s := \z(a) = \{ x \in \g \mid [x, a] = 
0\},
\ee
so that $\hat \S \subset \G$ is a closed Lie subgroup whose Lie 
algebra equals $\hat \s$. Observe that the restriction $\A^{-1}(a) \ra M$ 
is a principal $\hat \S$-bundle, hence we conclude that $\hat \T \subset \hat
\S$. Moreover, on $\hat M$, we have 
\[
\hat \al = 2 \kappa a
\] 
for some $\kappa \in \Om^1(\hat M)$ which by (\ref{eq:dA-hat}) satisfies 
$d\kappa = \pi^*(\om)$. In particular, the {\em Ambrose-Singer Holonomy 
theorem} implies that $\T_a = \exp(\F a) \subset \G$ is the identity 
component of $\hat \T$ which is thus a one dimensional (possibly non-regular) 
subgroup of $\hat \S$, and $\kappa$ yields the desired connection form on the 
principal $\hat \T$-bundle $\hat M \ra M$ which shows the first part.

Define $\cC_a \subset \cC$ as in (\ref{eq:cC0}) and $\Gamma_a \subset \G$ and 
$Q \subset \g$ as in (\ref{eq:Gamma0}), and let 
\be \label{eq:def-hat-B}
\hat B := p^{-1}(\hat M) \subset B \times \G,
\ee
where $p: B \times \G \ra B \times_\eH \G = \B$ is the canonical
projection. Then the restriction of the map 
\[
\ov \imath: B \times \G \longrightarrow \G,\ \ \ \ \ \ov \imath(b,g) := g^{-1}
\]
satisfies $\ov \imath(\hat B) \subset \Gamma_a$; indeed, since $\A(\hat M) 
\equiv a$, it follows that $\Ad_{g^{-1}} A(b) = a$ for all $(b,g) \in 
\hat B$ and hence $\Ad_g a = A(b) \in Q$, so that $g^{-1} \in \Gamma_a$. 
Since $2 \kappa a = \hat \al = \Ad_{g^{-1}} \sigma + \mu$, it follows by 
(\ref{eq:defsigma}) that
\[
\ov \imath^*(\mu) = -\Ad_g \mu = -2 \kappa \Ad_g a + \sigma = -2 \kappa A + e_- 
\ot \th + \eta + e_+ \ot (\rho \th) + \frac12 \om(u, \th) e_+^2,
\]
and hence
\[
\ov \imath^*(\mu) = -2 \kappa \left( \frac12 e_-^2 + \rho 
\right) + e_- \ot \th + \eta\ \ \ \ \mod \g^1 \oplus \g^2.
\]
Comparing this equation with the structure equations in
Theorem~\ref{thm:Maurer-Cartan}, it follows that the induced map $\hat 
\imath: \hat M = \hat B/\eH \ra \cC_a = \Gamma_a/\eH$ is a local 
diffeomorphism and the induced map $\imath: \tilde M := \T \backslash \hat M 
\ra \T_a \backslash \cC_a$ is connection preserving, where $\T_a \backslash 
\cC_a$ is (locally) equipped with the special symplectic connection from 
Theorem~\ref{thm:canonicalconn}.
{\hfill \rule{.5em}{1em}\mbox{}\bigskip}

\begin{rem} {\em
The proof of Theorem~B generalizes immediately to orbifolds. Namely, if $M$ is
an {\em orbifold}, then a special symplectic orbifold connection consists of 
an {\em almost principal $\eH$-bundle $B \ra M$}, i.e. $\eH$ acts locally 
freely and properly on $B$ such that $M = B/\eH$, and a coframing $\th + \eta 
\in \Om^1(B) \ot (V \oplus \h)$ on $B$ such that $\eta(\xi_h) \equiv h \in 
\h$ and $\th(\xi_h) \equiv 0$ for all infinitesimal generators $\xi_h$ of the 
$\eH$-action, and such that the structure equations (\ref{eq:structuresimple})
hold for some function $\rho: B \ra \h$.

Now the proofs of Propositions~\ref{prop:allareinduced} and 
\ref{prop:extendbundle} as well as the proof of Theorem~B go through verbatim 
as we never used the freeness of the $\eH$-action on $B$. In particular, the
holonomy reduction $\hat M$ is a {\em manifold} on which $\hat \T$ acts 
locally freely, and $M = \hat \T \backslash \hat M$ as an orbifold.}
\end{rem}

\section{Symmetries and compact special symplectic manifolds}

\begin{df} Let $(M, \nabla)$ be a manifold with a connection.
A {\em (local) symmetry} of the connection is a (local) diffeomorphism 
$\und \phi: M \ra M$ which preserves $\nabla$, i.e. such that 
$\nabla_{d\und \phi(X)}d\und \phi(Y) = d\und \phi(\nabla_X Y)$ for all vector 
fields $X, Y$ on $M$.
An {\em infinitesimal symmetry} of the connection is a vector field 
$\und \zeta$ on $M$ such that for all vector fields $X, Y$ on $M$ we have the 
relation
\[
{}[\und \zeta, \nabla_X Y] = \nabla_{[\und \zeta, X]} Y + \nabla_X 
[\und \zeta, Y].
\]

Furthermore, let $\pi: B \ra M$ be an $\eH$-structure compatible 
with $\nabla$, and let $\th, \eta$ denote the tautological and the connection 
form on $B$, respectively. A {\em (local) symmetry} on $B$ is a 
(local) diffeomorphism $\phi: B \ra B$ such that $\phi^*(\th) = \th$ 
and $\phi^*(\eta) = \eta$.
An {\em infinitesimal symmetry} on $B$ is a vector field $\zeta$ on 
$B$ such that ${\frak L}_\zeta (\th) = {\frak L}_\zeta (\eta) = 0$.
\end{df}

The ambiguity of the terminology above is justified by the one-to-one 
correspondence between (local or infinitesimal) symmetries on $M$ 
and $B$. Namely, if $\und \phi: M \ra M$ is a (local) symmetry, then 
there is a unique (local) symmetry $\phi: B \ra B$ with $\pi \circ \phi = 
\und \phi \circ \pi$, and vice versa. Likewise, for any infinitesimal 
symmetry $\und \zeta$ on $M$, there is a unique infinitesimal symmetry 
$\zeta$ on $B$ such that $\und \zeta = d\pi(\zeta)$.

The infinitesimal symmetries form the Lie algebra of the (local) group 
of (local) symmetries. We also observe that an infinitesimal symmetry 
on $B$ is uniquely determined by its value at any point. (The 
corresponding statement fails for infinitesimal symmetries on $M$ in 
general.)

\

\noindent
{\bf Proof of Corollary~C.}
The first part follows immediately from Theorem~B since $\cC_a \subset \cC$ 
is an open subset of the analytic manifold $\cC$, and the action of $\T_a$ on 
$\cC_a$ is analytic as well. Also, the $C^4$-germ of the connection at a 
point determines uniquely the $\G$-orbit of $a \in \g$ by 
(\ref{eq:structureGamma}) and hence the connection by Theorem~B.

Note that the generic element $a \in \g$ is $\G$-conjugate to an element in 
the Cartan subalgebra which is uniquely determined up to the action of the
(finite) Weyl group. Since multiplying $a \in \g$ by a scalar does not change
the connection, it follows that the generic special symplectic connection
associated to $\g$ depends on $(\rk(\g) -1)$ parameters.

For the second part, by virtue of Theorem~B it suffices to show the statement 
for manifolds of the form $M = M_U$ where $U \subset \cC_a$ is a regular 
open subset for some $a \in \g$. Let $\Gamma_U \subset \Gamma_a \subset \G$ 
be the $\eH$-invariant subset such that we have the principal $\eH$-bundle
$\Gamma_U \ra U$, and let $B_U:= \T_a \backslash \Gamma_U$ so that $B_U \ra 
M_U$ is the associated $\eH$-structure. 

Let $x \in \hat \s$, and denote by $\hat \zeta_x$ the right invariant vector 
field on $\G$ corresponding to $-x$, so that the map $x \mapsto \hat 
\zeta_x$ is a Lie algebra homomorphism. Then ${\frak L}_{{\hat \zeta}_x} 
(\mu) = 0$ where $\mu$ denotes the Maurer-Cartan form. By (\ref{eq:Gamma0}), 
it follows that the restriction of $\hat \zeta_x$ to $\Gamma_a$ is tangent, 
and since $\Gamma_U \subset \Gamma_a$ is open, we may regard $\hat \zeta_x$ 
as a vector field on $\Gamma_U$. Since $\hat \zeta_x$ commutes with the 
action of $\T_a$, it follows that there is a related vector field $\zeta_x$ 
on the quotient $B_U = \T_a^{loc} \backslash \Gamma_U$, and since the 
tautological and curvature form of the induced connection on $B_U$ pull back 
to components of $\mu$, it follows that $\zeta_x$ is an infinitesimal 
symmetry on $B_U$. 

Conversely, suppose that $\zeta$ is an infinitesimal symmetry on $B_U$. Since 
an infinitesimal symmetry must preserve the curvature and its covariant 
derivatives, we must have $\zeta(A) = 0$. But the tangent of the fiber of the 
map $A: \Gamma_U \ra \g$ is spanned by the vector fields $\zeta_x$, $x \in 
\hat \s$, and since infinitesimal symmetries are uniquely determined by their 
value at a point, it follows that $\zeta = \zeta_x$ for some $x \in \hat \s$.

Finally, it is evident that $\zeta_x = 0$ iff $\hat \zeta_x$ is tangent to 
$\T_a$ iff $x \in \F a$, hence the claim follows.
{\hfill \rule{.5em}{1em}\mbox{}\bigskip}

The rest of this section shall be devoted to the study of {\em compact simply
connected} manifolds with special symplectic connections. In fact, the main
result which we aim to prove is the following

\begin{thm} \label{thm:compact}
Let $\g$ be a $2$-gradable simple Lie algebra, let $\G$ be the connected Lie 
group with Lie algebra $\g$ and trivial center, and let $\hat \S \subset \G$ 
be a maximal compact subgroup. Then $\cC = \hat \S/\K$ for some compact 
subgroup $\K \subset \hat \S$ where $\cC \subset \P^o(\g)$ is the root cone. 
Moreover, let $\T \subset \hat \S$ be the identity component of the center of 
$\hat \S$. Then the following are equivalent:
\bi
\item
There is a compact simply connected symplectic manifold $M$ with a special
symplectic connection associated to the simple Lie algebra $\g$.
\item
$\g$ is a {\em real} Lie algebra and $\dim \T = 1$, i.e. $\T \cong S^1$.
\item
$\g$ is a {\em real} Lie algebra and $\T \neq \{e\}$.
\ei

If these conditions hold then $\T \backslash \cC \cong \hat \S/(\T \cdot 
\K)$ is a compact hermitian symmetric space, and the map $\imath: M \ra \T
\backslash \cC$ from Theorem~B is a connection preserving covering. Thus, $M$ 
is a hermitian symmetric space as well.
\end{thm}

This theorem allows us to classify all compact simply connected 
manifolds with special symplectic connections, as the maximal compact 
subgroups of semisimple Lie groups are fully classified (e.g. \cite{OV}). 
Thus, we obtain Theorem~D from the introduction as an immediate consequence.

The proof of Theorem~\ref{thm:compact} will be split up into several
steps. First, we observe the following

\begin{lem} \label{lem:transitive}
If the connected Lie group $\G$ acts transitively on the compact manifold $X$,
then so does any maximal compact subgroup $\hat \S \subset \G$.
\end{lem}

Thus, we can write the root cone as $\cC = \hat \S/\K$ for some compact 
subgroup $\K \subset \hat \S$ as asserted in Theorem~\ref{thm:compact}.

\begin{proof}
Let $X = \G/\eH$ as a homogeneous space, and let $\K \subset \eH$ be a maximal
compact Lie subgroup. Then there is a maximal compact Lie subgroup $\hat \S 
\subset \G$ which contains $\K$. Since the inclusions $\hat \S 
\hookrightarrow \G$ and $\K \hookrightarrow \eH$ are homotopy equivalences, 
standard homotopy arguments imply that the inclusion $\hat \S/\K 
\hookrightarrow X$ is also a homotopy equivalence. In particular, since both 
spaces are compact, they have equal dimension, so that $\hat \S/\K = \G/\eH$.
\end{proof}

Let us now suppose that $M$ is {\em real}. The proof that the first
condition in Theorem~\ref{thm:compact} implies the second and that in this 
case $M$ is the universal cover of the hermitian symmetric space $\T 
\backslash \cC$ is pursued in Lemmas~\ref{lem:T=S1} through
Proposition~\ref{prop:1->3}.

\begin{lem} \label{lem:T=S1}
Let $M$ be a compact real simply connected manifold with a special symplectic
connection associated to the real Lie algebra $\g$, and let $a \in \g$, $\T_a 
\subset \G$ and $\cC_a \subset \cC$ as in Theorem~\ref{thm:canonicalconn}.
Then $\T_a \cong S^1$ and $\cC_a = \cC$. Moreover, $\T_a$ acts freely on the 
universal cover $\tilde \cC$ of $\cC$, and $M = \T_a \backslash \tilde \cC$.
\end{lem}

\begin{proof}
If $M$ is simply connected, then by Theorem~B from the introduction
there is an $a \in \g$ and a principal $\T_a$-bundle $\pi: \hat M \ra M$ with 
a connection form $\kappa$ whose curvature equals $d\kappa = \pi^*(\om)$. 
Thus, $\pi^*(\om)$ is exact, while $\om$ cannot be exact if $M$ is compact. 
This implies that $\pi$ cannot be a homotopy equivalence, i.e. $\T_a$ cannot 
be contractible, hence $\T_a \cong S^1$. Thus, $\hat M$ is also compact, 
hence the local diffeomorphism $\hat \imath: \hat M \ra \cC$ from 
(\ref{diagram2}) must be surjective and a finite $\T_a$-equivariant covering. 
In particular, $\cC_a = \cC$.

Therefore, there is a $\T_a$-equivariant covering $\tilde \cC \ra \hat M$, 
where $\tilde \cC$ is the universal cover of $\cC$, and since $\T_a$ acts 
freely on $\hat M$, it acts also freely on $\tilde \cC$. Thus, the induced 
map $\T_a \backslash \tilde \cC \ra \T_a \backslash \hat M = M$ must also be 
a covering, hence a diffeomorphism as $M$ is simply connected.
\end{proof}

\noindent
We continue with the investigation of two special classes of examples.

\begin{prop} \label{prop:Bochner}
For $\g := \su(p+1, q+1)$ with $p+q \geq 1$, there are two orbits of maximal 
root vectors which are negatives of each other and are hence denoted by $\cC$ 
and $-\cC$. Moreover, these orbits are simply connected.

Let $a \in \g$ be such that $\T_a \cong S^1$, $\cC_a = \cC$ and 
the action of $\T_a$ on $\cC_a$ is free. Then $a$ is conjugate to a  
scalar multiple of $diag((q+1)i, \ldots, (q+1)i, -(p+1)i, \ldots, -(p+1)i)$. 
In particular, $\T_a \backslash \cC \cong \C\P^p \times \C\P^q$ with the 
hermitian symmetric connection as described in Theorem~E.
\end{prop}

\begin{proof}
We let $J: \C^{p+1,q+1} \ra \C^{p+1,q+1}$ be the $\g$-equivariant complex 
structure such that the metric $g(x,y) := \om(Jx, y)$ has signature $(p+1, 
q+1)$. Now $\g$ is a real form of $\sl(p+q+2, \C)$ whose maximal root cone 
consists of all traceless endomorphisms of (complex) rank $1$, hence the same is true 
for $\g$. The image of such an endomorphism must be a null line, so that the
maximal root cone of $\g$ consist of all endomorphisms of the form
\[
\{ \al_x \mid x \neq 0, g(x,x) = 0\} \dot \cup \{ -\al_x \mid x \neq 0, g(x,x) = 0\} 
=: \hat \cC \dot \cup (- \hat \cC),
\]
where
\[
\al_x(v) := g(v, x) Jx - g(v, Jx) x.
\]
Observe that $\al_{\la x} = |\la| \al_x$ for all 
$\la \in \C^*$, hence the projectivizations $\pm \cC$ of $\pm \hat \cC$ 
consist of all null lines in $\C^{p+1,q+1}$.

Decomposing $\C^{p+1,q+1} = \C^{p+1,0} \oplus \C^{0,q+1} =: \C^+ \oplus 
\C^-$, each null vector can be written as $x = x_+ + x_-$ with $x_\pm \in 
\C^\pm$ and $||x_+|| = ||x_-||$. In particular, $\cC = (S^{2p+1} \times
S^{2q+1})/diag(S^1)$, and a glance at the homotopy exact sequence now implies 
that $\cC$ is simply connected for $p+q \geq 1$.

Let $a \in \g$ be such that $\T_a \cong S^1$. Then $a$ is conjugate to an 
element of the form 
\linebreak
$diag(i \th_0, \ldots i \th_p, i \psi_0, \ldots, i 
\psi_q)$. If we denote the standard basis of $\C^{p+1, q+1}$ by $e_0,\ldots, 
e_p, f_0, \ldots, f_q$, then $x = e_r + f_s$ is a null vector and $(a, x 
\circ x) = \th_r - \psi_s$. Since $x \circ x \in \hat \cC$, this implies that 
$\th_r > \psi_s$ for all $r,s$.

Consider $T := \exp(2\pi/(\th_r - \psi_s) a) \in \T_a$. We have $T(e_r + 
f_s) = \exp(2\pi i\th_r/(\th_r - \psi_s)) (e_r + f_s)$ so that $T$ fixes $\C 
(e_r + f_s) \in \cC$. Thus, since $\T_a$ acts freely on $\cC$, it follows 
that $T = \exp(2\pi i\th_r/(\th_r - \psi_s)) Id$, which implies that 
$\exp(2\pi i\th_t/(\th_r - \psi_s)) = \exp(2\pi i\psi_u/(\th_r - \psi_s)) = 
\exp(2\pi i\th_r/(\th_r - \psi_s))$ for all $t,u$. Therefore, $(\th_t - 
\psi_u)/(\th_r - \psi_s) \in \Z$ for all $r,s,t,u$, and by switching $(r,s)$ 
and $(t,u)$ we conclude that $(\th_t - \psi_u)/(\th_r - \psi_s) = \pm 1$. But 
$\th_t - \psi_u, \th_r - \psi_s > 0$, whence this quotient must equal $1$ for 
all $r,s,t,u$, so that $\th_r = \th_t$ and $\psi_s = \psi_u$ for all 
$r,s,t,u$, hence $a$ must be of the asserted form, and the remaining 
statements now follow from the construction of the special symplectic 
connection.
\end{proof}

Evidently, there is no need to consider both maximal root orbits $\cC$ and $-\cC$, 
since one is obtained from the other by replacing the symplectic form $\om$ 
(or the complex structure $J$) by its negative which is irrelevant for our purposes. 
An analogous remark applies to the following case.

\begin{prop} \label{prop:RicciType}
For $\g := \sp(n+1, \R)$ with $n \geq 2$, there are two orbits of maximal 
root vectors which are negatives of each other and are hence denoted by $\cC$ and 
$-\cC$. Moreover, $\cC$ and $-\cC$ are diffeomorphic to $\R\P^{2n+1}$ and
therefore have fundamental group $\Z_2$.

Let $a \in \g$ be such that $\T_a \cong S^1$, $\cC_a = \cC$ and the action of $\T_a$ 
on the universal cover $\tilde \cC \cong S^{2n+1}$ is free. Then $a = c J$ for some 
$c > 0$, where $J$ is a complex structure on $\R^{2n+2}$ such that $g(x,y) := \om(Jx, y)$ is 
symmetric and positive definite. In particular, $\T_a \backslash \tilde \cC = 
\T_a \backslash \cC \cong \C\P^n$ with the hermitian symmetric connection as 
stated in Theorem~E.
\end{prop}

\begin{proof}
Since the conjugatets of the maximal root vectors in $\sp(n+1,\R)$ are the elements of rank one, 
the maximal root cone can be written as $\{ x \circ x \mid x \neq 0\} \dot \cup  
\{ -x \circ x \mid x \neq 0\} =: \cC \dot\cup -\cC$, where the product 
$\circ: S^2(\R^{n+1}) \ra \sp(n+1, \R)$ is given in (\ref{eq:circle-spV}). Thus, $\pm \cC = 
\pm \hat \cC/\R^+ \cong \R\P^{2n+1}$.

Let $J \in \sp(n+1, \R)$ be a complex structure such that $\om(Jx, x) > 0$ for
all $x \neq 0$, and let $a \in \g$ be such that $\T = \exp(\R a) \cong 
S^1$. Then $a$ is conjugate to an element of the form $J diag(\th_1, \ldots
\th_{2n+2})$, and $\cC_a = \cC$ implies that $0 < (a, x \circ x) = \om(a 
x, x)$ which is equivalent to $\th_i > 0$ for all $i$.

Consider $T := \exp(\pi/\th_i a) \in \T_a$. We have $T(e_i) = -e_i$ so that 
$T$ (if we consider the action on $\cC$) or $T^2$ (if we consider the action 
on $\tilde \cC$) has a fixed point. Thus, it follows that $T(e_j) = \pm e_j$ 
for all $j$ which implies that $\th_j | \th_i$ for all $j$, and switching the 
roles of $i$ and $j$, it follows that $\th_i = \pm \th_j$. But since $\th_i > 
0$ for all $i$, we must have $\th_i = \th_j$, hence $a$ is of the asserted 
form, and the remaining statements follow.
\end{proof}

\noindent
In order to work towards the general case, we continue with the following

\begin{lem} \label{lem:A0Conjugate}
Let $\g$ be a real $2$-gradable simple Lie algebra with the decomposition
(\ref{eq:decomposeg}). Let $a \in \g$ be such that $\T_a \cong S^1$ is a 
circle. Then $a$ is conjugate to an element of the form
\be \label{A0conjugate}
\frac c2 (e_+^2 + e_-^2) + \rho_0
\ee
with $c \in \R$ and $\rho_0 \in \h$.
\end{lem}

\begin{proof}
Let $\T_G \subset \G$ be the maximal compact abelian subgroup containing the 
circle $\SO(2) := \exp(\R(e_+^2 + e_-^2))$. Since $stab(e_+^2 + e_-^2) = 
\R(e_+^2 + e_-^2) \oplus \h$, it follows that $\T_G \subset \SO(2) \cdot 
\eH$.

The lemma now follows since {\em any} subgroup of $\G$ isomorphic to $S^1$ is
conjugate to a subgroup of $\T_G$.
\end{proof}

\begin{lem} \label{lem:g-beta}
Let $\g$ be a real $2$-gradable simple Lie algebra with the decomposition
(\ref{eq:decomposeg}), and let $\al_0 \in \Delta$ be the long root with
$\g^{\pm2} = \g_{\pm \al_0}$. Let $\beta \in \Delta$ be a root with 
$\pair \beta {\al_0} = 1$, and let $\ov \beta$ denote the conjugate root
w.r.t. the real form $\g$. If we define
\[
\g_{<\beta>} := \g \cap \left< \g_{\pm \al_0} \oplus \g_{\pm \beta} \oplus 
\g_{\pm \ov \beta} \right>,
\]
where $<\ >$ denotes the generated Lie subalgebra, then $\g_{<\beta>}$ is 
isomorphic to either $\sl(3,\R)$, $\sp(2,\R)$, $\g_2'$, $\su(1,2)$, or 
$\so(2,4) \cong \su(2,2)$.
\end{lem}

\begin{proof}
Since $\al_0$ is a real root, it follows that $\{\al_0, \beta, \ov \beta\}$ is
invariant under conjugation, hence $\g_{<\beta>}$ is a real form of the 
complex simple Lie algebra whose root system is generated by $\{\al_0, \beta,
\ov \beta\}$. Since $\g_{\pm \al_0} \subset \g_{<\beta>}$, it follows that
$\g_{<\beta>}$ is also $2$-gradable, and the decomposition
(\ref{eq:decomposeg}) reads $\g_{<\beta>} = \bigoplus_{i=-2}^2 (\g_{<\beta>} 
\cap \g^i)$.

If $\beta = \ov \beta$ is a real root, then the root system generated by 
$\al_0, \beta$ is irreducible of rank two and contains only real roots,
i.e. $\g_{<\beta>}$ is the split real form of type $A_2$, $B_2$ or $G_2$ as
listed above.

Therefore, for the rest of the proof we shall assume that $\beta \neq \ov 
\beta$. Since $\al_0$ is real, it follows that $\pair {\ov \beta}{\al_0} = 
\pair \beta {\al_0} = 1$, hence $\ov \beta \neq - \beta$ and thus $\beta, 
\ov \beta$ are linearly independent roots of equal length, so that they 
generate a root system either of type $A_2$ or of type $A_1 + A_1$. 
Since this root system is invariant under conjugation, it follows that there 
is a corresponding subalgebra $\hat \g_{<\beta>} \subset \g_{<\beta>}$ which 
is a real form of either $\sl(3,\C)$ or $\so(4,\C)$. This real form must 
contain roots which are neither real nor purely imaginary since $\beta \neq 
\pm \ov \beta$. In particular, it is neither split nor compact, and thus, the 
only real forms possible are $\su(1,2)$ in the first and $\so(1,3)$ in the 
second case.

If $\hat \g_{<\beta>} \cong \su(1,2)$, then $\hat \g_{<\beta>}$ is 
$2$-gradable, and hence the root system generated by $\beta, \ov \beta$ must 
contain a real root. This implies that $\beta + \ov \beta \in \Delta$, and 
since $\pair {\beta + \ov \beta}{\al_0} = 2$, it follows that $\beta + \ov 
\beta = \al_0$, i.e. $\hat \g_{<\beta>} = \g_{<\beta>} \cong \su(1,2)$.

Let us now suppose that $\hat \g_{<\beta>} \cong \so(1,3)$. We assert that in
this case, $\beta$ must be a long root. For if $\beta$ and hence $\ov \beta$ 
are short, then $\pair \beta {\al_0} = \pair {\ov \beta}{\al_0} = 1$ implies 
that $\beta + \ov \beta = \al_0$ so that the root system generated by 
$\{\al_0, \beta, \ov \beta\}$ is irreducible of rank two with roots of 
different length, i.e. $\g_{<\beta>}$ is a $2$-gradable real form with 
root system $B_2$ or $G_2$. However, by Table~1, the only $2$-gradable real 
forms of these root systems are the split forms which have only real 
roots, contradicting that $\beta \neq \ov \beta$.

Thus, we are left with the case where $\hat \g_{<\beta>} \cong \so(1,3)$ and 
$\beta \in \Delta$ is a long root. Then $\pair \beta {\ov \beta} = 0$, and 
the intersections 
\[
W_\pm := \g \cap (\g_{\pm \al_0} \oplus \g_{\pm(\al_0 - \beta) } \oplus 
\g_{\pm(\al_0 - \ov \beta)} \oplus \g_{\pm(\al_0 - \beta - \ov \beta)})
\]
are $\hat \g_{<\beta>}$-modules. In fact, considering the weights of the 
action of $\hat \g_{<\beta>}$ on $W_\pm$ implies that $W_\pm \cong \R^{1,3}$ 
as a $\hat \g_{<\beta>}$-module, and one verifies that $[W_+, W_+] = [W_-, 
W_-] = 0$, whereas $[W_+, W_-] \subset \hat \g_{<\beta>} \oplus \R 
H_{\al_0}$. It follows now that $(\g_{<\beta>}, \hat \g_{<\beta>} \oplus \R 
\eH_{\al_0})$ is an irreducible symmetric pair whose isotropy representation
coincides with that of the symmetric pair $(\so(2,4), \so(1,3) \oplus 
\so(1,1))$, hence these symmetric pairs are isomorphic. In particular, we 
have $\g_{<\beta>} \cong \so(2,4) \cong \su(2,2)$ which completes the proof.
\end{proof}

\begin{lem} \label{lem:G-beta}
Let $\g$ be one of the real Lie algebras from Lemma~\ref{lem:g-beta}, and let 
$a \in \g$ be such that $\T_a \cong S^1$. Define $\cC_a \subset \cC$ as in 
(\ref{eq:cC0}). Then
\bi
\item
If $\g \cong \sl(3,\R)$, $\g_2'$ then $\cC_a \subsetneq \cC$ is a {\em 
proper} subset for any such $a \in \g$.
\item
Let $\g \cong \sp(2,\R)$, $\su(1,2)$, $\su(2,2)$ and $\tilde \cC$ be the 
universal cover of $\cC$. If $\cC_a = \cC$ and the (lifted) action of $\T_a$ 
on $\tilde \cC$ is free, then the action of $\T_a$ on $\cC$ is free and $a$ 
is conjugate to an element of the form (\ref{A0conjugate}) with $c > 0$, 
$\rho_0^2 = -c^2 Id_V$ and $\om(\rho_0 x, x) > 0$ for all $0 \neq x \in V$.
\ei
\end{lem}

\begin{proof}
By Lemma~\ref{lem:A0Conjugate}, we may assume that $a$ is of the form 
(\ref{A0conjugate}). Since $\sl(3,\R)$, $\g_2'$ are split real forms, it 
follows that $\g_\beta \subset V_1 \cap \cC$ for all long roots $\beta$ with 
$\pair \beta {\al_0} = 1$, hence $V_1 \cap \cC \neq \emptyset$, whereas 
$(a, V_1 \cap \cC) = 0$. Thus, $\cC_a \neq \cC$.

The second part now follows immediately from
Propositions~\ref{prop:Bochner} and \ref{prop:RicciType} where the explicit
form of $a$ was given.
\end{proof}

\noindent
This lemma now allows us to treat the general case. Namely we have

\begin{lem} \label{lem:a-complex structure}
Let $\g$ be a $2$-gradable real Lie algebra, and let $a \in \g$ be 
such that $\T_a \cong S^1$, $\cC_a = \cC$ and that the action of $\T_a$ on 
the universal cover of $\cC$ is free. Then $a$ is conjugate to an element of 
the form (\ref{A0conjugate}) with $c > 0$, $\rho_0^2 = -c^2 Id_V$ and 
$\om(\rho_0 x, x) > 0$ for all $0 \neq x \in V$.
\end{lem}

\begin{proof}
Lemma~\ref{lem:A0Conjugate} allows us to assume that $a$ is of the form 
(\ref{A0conjugate}). Indeed, we may assume that $\rho_0$ is contained in the 
Cartan subalgebra of $\h_\C$, so that the Lie subalgebras $\g_{<\beta>} 
\subset \g$ from Lemma~\ref{lem:g-beta} are $\T_a$-invariant.

Let $\G$ be a connected Lie group with Lie algebra $\g$ and let $\G_{<\beta>} 
\subset \G$ be the connected Lie subgroup with Lie algebra $\g_{<\beta>} 
\subset \g$. Then $\cC_\beta := \G_{<\beta>} \cdot e_+^2 \subset \cC$ is 
$\T_a$-invariant, and $(\cC_\beta)_a = \cC_\beta \cap \cC_a = \cC_\beta$ as 
$\cC_a = \cC$. Thus, since $\cC_\beta$ is the cone of maximal roots of
$\g_{<\beta>}$, by Lemma~\ref{lem:g-beta} and the first part of 
Lemma~\ref{lem:G-beta} we conclude that $\g_{<\beta>} \cong \sp(2,\R), 
\su(1,2)$ or $\su(2,2)$.

The inverse image $\hat M_{<\beta>} := \hat \imath^{-1}(\cC_\beta) \subset 
\hat M$ with the covering $\hat \imath: \hat M \ra \cC$ from (\ref{diagram2})
must also be $\T_a$-invariant, and every connected component of 
$\hat M_{<\beta>}$ is a covering of $\cC_\beta$. Since $\T_a$ acts freely on 
$\hat M_{<\beta>} \subset \hat M$, it follows from the second part of
Lemma~\ref{lem:G-beta} that $c > 0$, 
$\rho_0^2|_{V_{<\beta>}} = -c^2 Id$ and $\om(\rho_0 x, x) > 0$ for all $0 
\neq x \in V_{<\beta>}$, where $V_{<\beta>} \subset V$ is defined by the 
relation $e_\pm \ot V_{<\beta>} = \g_{<\beta>} \cap \g^{\pm1}$.

The claim now follows since $V$ is the direct sum of the $V_{<\beta>}$, and 
for all $\beta, \gamma \in \Delta$ with $\pair \beta {\al_0} = \pair
\gamma {\al_0} = 1$ and $V_{<\beta>} \cap V_{<\gamma>} = 0$ we have
$\om(V_{<\beta>}, V_{<\gamma>}) = 0$.
\end{proof}

\begin{lem} \label{lem:hat s}
Let $\g$ be a $2$-gradable real Lie algebra, and let $a \in \g$ be of 
the form (\ref{A0conjugate}) with $c > 0$, $\rho_0^2 = -c^2 Id_V$ and 
$\om(\rho_0 x, x) > 0$ for all $0 \neq x \in V$. Then 
\be \label{eq:hat s}
\hat \s = stab(a) := \{ x \in \g \mid [x, a] = 0\} = \R a \oplus 
\k \oplus \{ c (e_+ \ot x) - (e_- \ot \rho_0 x) \mid x \in V\},
\ee
where $\k := \{ h \in \h \mid [h, \rho_0] = 0\}$. Moreover, $\hat \s \cong \R 
a \oplus \s$, where $\s$ is a compact semisimple Lie algebra, and $(\hat \s, 
\R a \oplus \k)$ is a hermitian symmetric pair. Also, $\hat \s \subset \g$ is 
a maximal Lie subalgebra.
\end{lem}

\begin{proof}
It is straightforward to verify (\ref{eq:hat s}) and $\z(\hat \s) = \R a$, and 
that $(\hat \s, \R a \oplus \k)$ is a hermitian symmetric pair. Also, note 
that $\k$ is the Lie algebra of the compact group $\K = \eH \cap \U(V, 1/c\ 
\rho_0)$. Thus, there is a positive definite $\ad_\k$-invariant metric on 
$\g$, so that $\ad_h^2: \g \ra \g$ is negative semidefinite for all $h \in \k$ 
and hence $B(h,h) = tr(\ad_h^2) \leq 0$ with equality iff $\ad_h = 0$ iff $h = 
0$ since $\g$ is simple and hence has trivial center. Thus, $(h,h) > 0$ for 
all $0 \neq h \in \k$ by (\ref{eq:Killing}). Also, $((e_+ \ot x) - (e_- \ot 
\rho_0 x), (e_+ \ot x) - (e_- \ot \rho_0 x)) = 2 \om(\rho_0 x, x) > 0$ for all $0 
\neq x \in V$, and $(e_+^2 + e_-^2, e_+^2 + e_-^2) = 4 > 0$. Since $e_+^2 + 
e_-^2$, $\k$ and $\{c (e_+ \ot x) - (e_- \ot \rho_0 x) \mid x \in V\}$ are 
orthogonal w.r.t. $(\ ,\ )$, it follows that $(\ ,\ )$ is positive definite and 
$\ad$-invariant on $\hat \s$. Thus, $\ad_x: \hat \s \ra \hat \s$ is skew 
symmetric w.r.t. $(\ ,\ )$ for all $x \in \hat \s$, so that $B_{\hat \s}(x,x) = tr(\ad_x^2) 
\leq 0$ with equality iff $x \in \z(\hat \s)$, hence $\hat \s = \z(\hat \s) \oplus \s$ 
for a compact semisimple Lie algebra $\s$ as asserted.

To see that $\hat \s \subset \g$ is a maximal subalgebra, let $\hat \s 
\subset \g' \subsetneq \g$ be a subalgebra. Considering the eigenspaces of 
$\ad(e_+^2 + e_-)^2$, it follows that $\g' = (\g' \cap \sl(2,\R)) \oplus (\g' 
\cap \h) \oplus (\g' \cap \R^2 \ot V)$.

But $\sl(2,\R)$ and $\hat \s$ generate $\g$, so it follows that $\g' \cap 
\sl(2,\R) = \R(e_+^2 + e_-^2)$. Also, if $e_+\ot x \in \g'$, then $[e_+ 
\ot x, e_+ \ot y - e_- \ot \rho_0 y] \in \g'$ implies $\om(x,y) = 0$, as one 
sees by looking at the $\sl(2,\R)$-component. Since this is the case for all 
$y \in V$, it follows that $\g' \cap \R^2 \ot V = \hat \s \cap \R^2 \ot V$. 
Finally, if $h \in \g' \cap \h$ then $[h, e_+ \ot x - e_- \ot \rho_0 x] \in 
\g' \cap \R^2 \ot V \subset \hat \s$, and from here it follows that $h \in 
\k$, so that $\g' = \hat \s$ as claimed.
\end{proof}

Now we are ready to prove that in the real case, the first condition in 
Theorem~\ref{thm:compact} implies the second, and that in this case
$M$ is hermitian symmetric.

\begin{prop} \label{prop:1->3} 
Let $M$ be a real compact simply connected manifold with a special symplectic
connection, and let $a \in \g$ be from Theorem~B. Then $\T_a \backslash \cC$ 
is a hermitian symmetric space, and the map $\imath: M \ra \T_a \backslash 
\cC$ is a connection preserving covering. Moreover, $\T_a$ is the connected
component of the center of $\hat \S \subset \G$, where $\hat \S$ is a maximal
compact subgroup of $\G$.
\end{prop}

\begin{proof}
By Lemma~\ref{lem:hat s}, it follows that the connected Lie subgroup $\hat \S 
\subset \G$ with Lie subalgebra $\hat \s$ must be compact as $\T_a \cong S^1$
is compact. Indeed, it is a maximal compact subgroup as $\hat \s \subset \g$ 
is maximal, and $\T_a \subset \hat \S$ is the connected component of its 
center.

Thus, if we write $\cC = \hat \S/\K$ by Lemma~\ref{lem:transitive}, then $\K$
has $\k = \hat \s \cap \p$ as its Lie algebra by (\ref{eq:hat s}), and hence 
$\T_a \backslash \cC = \hat \S/(\T_a \cdot \K)$ is a hermitian symmetric space
by Lemma~\ref{lem:hat s}, and the covering $\imath: M \ra \T_a \backslash 
\cC$ is connection preserving.
\end{proof}

Evidently, the second condition in Theorem~\ref{thm:compact} implies the 
third, hence the real case will be finished with the following 

\begin{lem} \label{lem:3implies1}
Let $\G$ be a real simple connected Lie group with $2$-gradable Lie algebra 
$\g$ and trivial center, and let $\hat \S \subset \G$ be a maximal compact 
Lie subgroup whose center contains $\T_a = \exp(\R a)$, some $0 \neq a \in 
\g$. Then - after changing $a$ to its negative if necessary - we have $\cC_a 
= \cC$, and the action of $\T_a$ on $\cC$ is free. Moreover, $\T_a \backslash 
\cC$ has finite fundamental group.
\end{lem}

By Theorem~\ref{thm:canonicalconn}, it then follows that $\T_a \backslash 
\cC$ carries a special symplectic connection associated to $\g$, hence so 
does its universal cover $M := (\T_a \backslash \cC)^{\tilde{\mbox{ }}}$. 
Since $\T_a \backslash \cC$ is compact and has finite fundamental group, $M$ 
is compact as well. Thus, the lemma shows that the third condition in 
Theorem~\ref{thm:compact} implies the first.

\begin{proof}
Since $\G$ acts transitively on $\cC$, so does $\hat \S$ by 
Lemma~\ref{lem:transitive}, hence we can write $\cC = \hat \S/\K$ for some 
compact subgroup $\K \subset \hat \S$. Let $a \in \g$ be such that $\T = 
\T_a$ and consider the corresponding contact symmetry $a^*$ from 
(\ref{eq:xi0}). We assert that $a^*$ is {\em transversal}. For if there 
is a $p \in \cC$ with $(a^*)_p \in {\cal D}_p$, then $dL_g((a^*)_p) \in 
dL_g({\cal D}_p) = {\cal D}_{g \cdot p}$ for all $g \in \hat \S$ by the 
$\hat \S$-equivariance of the contact structure. On the other hand,
\[
dL_g((a^*)_p) = \left. \frac d{dt} \right|_{t=0} g \cdot \exp(t a)
\cdot p = \left. \frac d{dt} \right|_{t=0} \exp(t \Ad_g (a)) \cdot g \cdot 
p = (a^*)_{g \cdot p},
\]
since $\Ad_g (a) = a$ for $g \in \hat \S$. Thus, $(a^*)_{g \cdot p} \in
{\cal D}_{g \cdot p}$, and since $\hat \S$ acts transitively on $\cC$, it 
follows that $(a^*)_q \in {\cal D}_q$ for {\em all} $q \in \cC$. But 
$a^*$ is a contact symmetry, hence this implies that $a^* \equiv 0$ which 
is a contradiction.

Thus, $\la(a^*) \neq 0$ for all $\la \in \hat \cC$ and - after replacing 
$a$ by its negative if necessary - we may assume that $\la(a^*) > 0$ for 
all $\la \in \hat \cC$, so that $\cC_a = \cC$. Since $\T_a$ lies in the 
center of $\hat \S$ and $\G$ acts effectively on $\cC = \hat \S/\K$ as $\G$ 
has trivial center, it follows that $\T_a$ acts freely on $\cC$.

It now follows from Lemmas~\ref{lem:a-complex structure} and \ref{lem:hat s}
that $\T_a \subset \hat \S$ is the connected component of the center, hence 
the inclusion $\T_a \cdot \K \hookrightarrow \hat \S$ induces a map with 
finite cokernel between the fundamental groups. Now the homotopy exact 
sequence of the fibration $\T_a \cdot \K \hookrightarrow \hat \S \ra \hat 
\S/(\T_a \cdot \K) = \T_a \backslash \cC$ implies that $\T_a \backslash \cC$ 
has finite fundamental group as claimed. 
\end{proof}

\noindent
Finally, we need to deal with the {\em complex} case which we do in the 
following

\begin{prop}
There are no compact simply connected complex manifolds $M$ with a special 
symplectic connection associated to a complex simple Lie algebra $\g$.
\end{prop}

\begin{proof}
If $M$ is such a manifold, then as in the proof of Lemma~\ref{lem:T=S1}, we
conclude that the fibration $\hat M \ra M$ cannot be a homotopy equivalence, 
so that $\T_a \cong \C^*$ and hence $a \in \g$ is semisimple. This means that
the eigenvalues of $\ad_{a}$ are all linearly dependent over $\Q$, so that 
-- after replacing $a$ by a suitable non-zero multiple -- we may assume 
that all these eigenvalues are integers. Thus, we may choose the split real 
form $\g_\R \subset \g$ such that $a \in \g_\R$ and $\T_a^\R := \exp(\R a) 
\subset \g_\R$ is isomorphic to $\R$.

Let $\cC^\R \subset \P^0(\g_\R)$ be the projectivization of the root cone of 
$\g_\R$, and consider the {\em Hopf fibration} $pr: \P^0(\g_\R) \ra \P^0(\g)$ 
which maps each real line to the corresponding complex one. Then $pr(\cC^\R) 
\subset \cC$ is a regular submanifold which is diffeomorphic to either 
$\cC^\R$ or $\cC^\R/\Z_2$. In particular, the restriction $pr: \cC^\R \ra 
pr(\cC^\R)$ is a regular covering. Also, as the distribution ${\cal D}$ 
consists of {\em complex} subspaces, it follows that 
\[
pr\left(\cC^\R_a \right) = pr(\cC^\R) \cap \cC_a,
\]
so that the restriction $pr: \cC^\R_a \ra pr(\cC^\R) \cap \cC_a$ is also a
regular covering. In particular, $pr(\cC^\R_a) \subset \cC_a$ is a regular 
closed submanifold.

Recall the covering map $\hat \imath: \hat M \ra \cC_a$ from 
(\ref{diagram2}). Standard homotopy arguments show that there is a manifold 
$\hat M^\R$ and regular coverings $\hat \imath^\R: \hat M^\R \ra \cC_a^\R$ 
and $\tilde{pr}: \hat M^\R \ra \hat \imath^{-1}(pr(\cC^\R) \cap \cC_a)$ where
$\hat \imath^\R$ is equivariant w.r.t. the action of $\T_a^\R \subset \T_a$. 
Note that $\T_a^\R$ acts freely and properly discontinuously on $\hat M$ and 
hence also on $\hat M^\R$, so that $M^\R := \T_a^\R \backslash \hat M^\R$ is 
a manifold. Hence, we obtain the following commutative diagram, where the 
dotted lines indicate immersions which are regular covers of their images 
with a deck group of order at most $2$:

\[
\xymatrix{
\hat M^\R \ar[dd]^{\mbox{\footnotesize $\T_a^\R$}} \ar@{.>}[drr]_{\tilde 
{pr}} \ar[r]^{\hat \imath^\R} 
& \cC_a^\R \ar[dd]^{\mbox{\footnotesize $\T_a^\R$}} \ar@{.>}[drr]_{pr} 
\\
& & \hat M \ar[dd]^{\mbox{\footnotesize $\T_a$}} \ar[r]^{\hat \imath} & 
\cC_a \ar[dd]^{\mbox{\footnotesize $\T_a$}}
\\
M^\R \ar[r]^{\imath^\R} \ar@{.>}[drr] &
\T_a^\R \backslash \cC_a^\R \ar@{.>}[drr]
\\
& & M \ar[r]^{\imath} & \T_a \backslash \cC_a
}
\]

Thus, Theorem~B implies that $M^\R$ carries a special symplectic connection 
associated to $\g_\R$, and the principal $\T_a^\R$-bundle $\hat M^\R \ra M^\R$
coincides with the one given in that theorem. 

But the image of the covering $M^\R \ra M$ is a closed submanifold, and since 
we assume that $M$ is compact, it follows that $M^\R$ is compact. Thus, as in 
the proof of Lemma~\ref{lem:T=S1}, we conclude that $\T_a^\R \cong S^1$, which
is a contradiction as $\T_a^\R \cong \R$ by our choice of $\g_\R$.
\end{proof}



\begin{thebibliography}{99}

{\small
\bibitem[BC1]{Michel} {\sc F.Bourgeois, M.Cahen}, {\em A variational principle
    for symplectic connections}, J.Geom.Phys {\bf 30}, 233-265 (1999)

\bibitem[BC2]{Michel4} {\sc F.Baguis, M.Cahen}, {\em A construction of
    symplectic connections through reduction}, Let.Math.Phys {\bf 57}, 
    149-160 (2001)

\bibitem[Bo]{Bo} {\sc S.Bochner}, {\em Curvature and Betti numbers, II}, 
    Ann.Math {\bf 50} 77-93 (1949)

\bibitem[BR]{BR} {\sc N.R.O'Brian, J.Rawnsley}, {\em Twistor spaces},
  Ann.Glob.Anal.Geom {\bf 3}, 29 - 58 (1985)

\bibitem[Br1]{Br1} {\sc R. Bryant}, {\em Two exotic holonomies in dimension
    four, path geometries, and twistor theory},  Proc. Symp. in Pure Math.
    {\bf 53}, 33-88 (1991)

\bibitem[Br2]{Bochner} {\sc R.Bryant}, {\em Bochner-K\"ahler metrics}, 
    J.AMS {\bf 14} No.3, 623-715 (2001)

\bibitem[CGR]{Michel2} {\sc M.Cahen, S.Gutt, J.Rawnsley}, {\em 
    Symmetric symplectic spaces with Ricci-type curvature},
    G.Dito, D.Sternheimer (ed.), Conf\'erence Mosh\'e Flato 1999, Vol.II, 
    Math.Phys.Stud. 22, 81-91 (2000)

\bibitem[CGHR]{Michel3} {\sc M.Cahen, S.Gutt, J.Horowitz, J.Rawnsley},
    {\em Homogeneous symplectic manifolds with Ricci-type curvature},
    J.Geom.Phys. {\bf 38} 140-151 (2001)

\bibitem[CGS]{CGS} {\sc M.Cahen, S.Gutt, L.J.Schwachh\"ofer}, {\em 
Construction of Ricci-type connections by reduction and induction},
preprint, arXiv:math.DG/0310375

\bibitem[CMS]{CMS} {\sc Q.-S. Chi, S.A. Merkulov, L.J. Schwachh{\"o}fer}, 
    {\em On the Existence of Infinite Series of Exotic Holonomies}, Inv. 
    Math. {\bf 126}, 391-411 (1996)

\bibitem[FH]{FH} {\sc W. Fulton, J. Harris}, {\em Representation Theory},
    Graduate Texts in Mathematics 129, Springer-Verlag, Berlin, New-York 
    (1996)

\bibitem[Hu]{Hu} {\sc J.E. Humphreys}, {\em Introduction to Lie Algebras 
    and Representation Theory}, Springer-Verlag, Berlin, New York (1987)

\bibitem[K]{K} {\sc Y.Kamishima} {\em Uniformization of K\"ahler manifolds 
    with vanishing Bochner tensor}, Acta Math. {\bf 172} No.2, 299-308 (1994)

\bibitem[OV]{OV} {\sc A.L. Onishchik, E.B. Vinberg}, {\em Lie groups and Lie 
    Algebras}, Vol. 3, Springer-Verlag, Berlin, New York (1996)

\bibitem[MS]{MS} {\sc S.A. Merkulov, L.J. Schwachh\"ofer}, {\em 
    Classification of irreducible holonomies of torsion free affine 
    connections}, Ann.Math. {\bf 150}, 77-149 (1999); {\em Addendum: 
    Classification of irreducible holonomies of torsion-free affine 
    connections}, Ann.Math. {\bf 150}, 1177-1179 (1999)

\bibitem[S1]{Habil} {\sc L.J. Schwachh\"ofer}, {\em On the classification of
    holonomy representations}, Habilitationsschrift, Universit\"at Leipzig 
    (1998)

\bibitem[S2]{S} {\sc L.J. Schwachh\"ofer}, {\em Homogeneous connections with 
    special symplectic holonomy}, Math.Zeit. {\bf 238}, 655 -- 688 (2001)

\bibitem[S3]{Advances} {\sc L.J. Schwachh\"ofer}, {\em Connections with 
    irreducible holonomy representations}, Adv.Math. {\bf 160}, 1 -- 80 
    (2001)

\bibitem[V]{V2} {\sc I.Vaisman}, {\em Variations on the theme of twistor 
    spaces}, Balkan J.Geom.Appl. {\bf 3} 135 - 156 (1998)
}
\end{thebibliography}
\end{document}